\documentclass[aos,preprint]{imsart}

\RequirePackage[OT1]{fontenc}
\RequirePackage{amsthm,amsmath,amsfonts}
\RequirePackage[numbers]{natbib}
\RequirePackage[colorlinks,citecolor=blue,urlcolor=blue]{hyperref}
\RequirePackage{hypernat}
\RequirePackage[pdftex]{graphics}
\RequirePackage{graphicx}
\RequirePackage{float}

\startlocaldefs
\numberwithin{equation}{section}
\theoremstyle{plain}
\newtheorem{theorem}{Theorem}[section]
\newtheorem{lemma}{Lemma}[section]

\newtheorem{example}{Example}[section]
\newtheorem{remark}{Remark}[section]
\newtheorem{assumption}{Assumption}[section]

\startlocaldefs
\endlocaldefs

\begin{document}

\begin{frontmatter}

\title{A General Framework for Consistency of Principal Component Analysis}
\runtitle{General Framework for PCA Consistency}


\begin{aug}
\author{\fnms{Dan} \snm{Shen}\thanksref{m1,t1}\ead[label=e1]{dshen@email.unc.edu}},
\author{\fnms{Haipeng} \snm{Shen}\thanksref{t2}\ead[label=e2]{haipeng@email.unc.edu}}
\and
\author{\fnms{J. S.} \snm{Marron}\thanksref{t3}
\ead[label=e3]{marron@email.unc.edu}}

\thankstext{m1}{Corresponding Author}
\thankstext{t1}{Partially supported by NSF grant DMS-0854908}
\thankstext{t2}{Partially supported by NSF grants DMS-1106912 and CMMI-0800575, NIH Challenge Grant 1 RC1 DA029425-01, and the Xerox
Foundation UAC Award}
\thankstext{t3}{Partially supported by NSF grants DMS-0606577 and DMS-0854908}
\runauthor{Dan Shen, Haipeng Shen and J. S. Marron}

\affiliation{University of North Carolina at Chapel Hill}

\address{Department of Statistics and Operations Research\\
University of North Carolina at Chapel Hill\\
Chapel Hill, NC 27599\\
\printead{e1}\\
\phantom{E-mail:\ }\printead*{e2}\\
\phantom{E-mail:\ }\printead*{e3}
}
\end{aug}

\begin{abstract}
A general asymptotic framework is developed for studying consistency properties of principal component analysis (PCA).
Our framework includes several previously studied domains of asymptotics as special cases and allows one to investigate interesting connections and transitions among the various domains. More importantly, it enables us to investigate asymptotic scenarios that have not been considered before, and gain new insights into the consistency, subspace consistency and strong inconsistency regions of PCA and the boundaries among them. We also establish the corresponding convergence rate within each region. Under general spike covariance models, the dimension (or the number of variables) discourages the consistency of PCA, while the sample size and spike information
(the relative size of the population eigenvalues)
encourages PCA consistency.
Our framework nicely illustrates the relationship among
these three types of information in terms of dimension, sample size and spike size,
and rigorously characterizes how their relationships affect PCA consistency.
\end{abstract}

\begin{keyword}[class=AMS]
\kwd[Primary ]{62H25}
\kwd[; Secondary ]{62F12}
\end{keyword}

\begin{keyword}
\kwd{PCA}
\kwd{High Dimension}
\kwd{Spike Model}
\kwd{Consistency}
\end{keyword}

\end{frontmatter}


\section{Introduction}\label{sec:01}

Principal Component Analysis (PCA) is  an important visualization and dimension
 reduction tool which finds orthogonal directions reflecting maximal variation
in the data. This allows the low dimensional representation of data, by projecting data onto
these directions. PCA is usually obtained by an eigen decomposition of the sample variance-covariance matrix of the data. Properties of the sample eigenvalues
and eigenvectors have been analyzed under several domains of asymptotics.

In this paper, we develop a {\it general asymptotic framework} to explore interesting transitions among the various asymptotic domains. The general framework includes the traditional asymptotic setups as special cases, which allows careful study of the connections among the various setups, and more importantly it investigates scenarios that have not been considered before, and offers new insights into the \emph{consistency}
({{in the sense that the angle between estimated and population eigen direction tends to 0, or the inner product tends to 1}}) and \emph{strong-inconsistency}
({{where the angle tends to $\frac{\pi}{2}$, i.e., the inner product tends to 0}}) properties of PCA, along with some technically challenging convergence rates.

Existing asymptotic studies of PCA roughly fall into three domains:
\begin{enumerate}
\item[(a)] the {\bf classical}  domain of asymptotics, under which the sample size $n\rightarrow \infty$ and the dimension $d$ is fixed
(hence the ratio $\frac{n}{d}\rightarrow \infty$).
For example, see \cite{girshick1939sampling,lawley1956tests,anderson1963asymptotic,anderson1984introduction, jackson1991user}.
\item[(b)] the {\bf random matrix} theory domain, where both the sample size $n$ and the dimension $d$ increase to infinity, with the ratio $\frac{n}{d}\rightarrow c$, a constant mostly assumed to be within $(0, \infty)$.  Representative work includes
    \cite{biehl1994statistical,watkin1994optimal,reimann1996gaussian,hoyle2003pca}
 from the statistical physics literature, as well as
\cite{johnstone2001distribution,baik2005phase,baik2006eigenvalues,onatski2006asymptotic,paul2007asymptotics,nadler2008finite,johnstone2009consistency,
lee2010convergence,benaych2011eigenvalues}
 from the statistics literature.
\item[(c)] the {\bf high dimension low sample size (HDLSS)}  domain of asymptotics, which is based on the limit, as the dimension $d\rightarrow\infty$, with the sample size $n$ being fixed (hence the ratio $\frac{n}{d}\rightarrow 0$). HDLSS asymptotics was originally studied by \cite{casella1982limit},
and recently rediscovered by
\cite{hall2005geometric}.
PCA has been studied using the HDLSS asymptotics by
\cite{ahn2007high,jung2009pca}.%
\end{enumerate}


PCA consistency and (strong) inconsistency, defined in terms of angles,
 are important properties that have been studied before.
 A common technical device is the spike covariance model, initially introduced by Johnstone~\cite{johnstone2001distribution}.
 This model has been used in this context by, for example, Nadler~\cite{nadler2008finite}, Johnstone and Lu~\cite{johnstone2009consistency},
 and Jung and Marron~\cite{jung2009pca}. An interesting, more general model has been considered by
 Benaych-Georges and Nadakuditi~\cite{benaych2011eigenvalues}.

Under the spike model, the first few eigenvalues are much larger than the others. A {\it major point of the present paper} is that there are three critical features whose relationships drive the consistency properties of PCA, namely
\begin{enumerate}
\item[(1)] the {\it sample information}: the sample size $n$, which has a {\it positive} contribution to, i.e. {\it encourages}, the consistency of the sample eigenvectors.
\item[(2)] the {\it variable information}: the dimension $d$, which has a {\it negative} contribution to, i.e. {\it discourages}, the
consistency of the sample eigenvectors.
\item[(3)] the {\it spike information}: the relative sizes of the several leading eigenvalues, 
which also has a positive contribution to the consistency. 
\end{enumerate}

Our general framework considers increasing sample size $n$,
increasing dimension $d$, and increasing spike
information. {{It clearly characterizes how their relationships determine the regions of
 consistency and strong-inconsistency of PCA,
 along with the boundary in-between.
 In addition, our theorems demonstrate the transitions among the existing domains of asymptotics, and for the first time to the best of our knowledge,}} enable one to understand the connections among them. Note that the classical domain ((a) above) assumes increasing sample size $n$ while fixing dimension $d$; the random matrix domain ((b) above) assumes increasing sample size $n$ and increasing dimension $d$,
while fixing the spike information; the HDLSS domain ((c) above) fixes the sample size, and increases the dimension and the spike information; {{thus each of these three domains is a boundary case of our framework.}} Finally, our theorems also contain novel results on rates of convergence.

Sections~\ref{sec:02} and~\ref{multiple-sike-models} formally state very general theorems for the single and multiple component spike models, respectively. For {illustration purposes only}, in this section we first consider Examples~\ref{example:01} and~\ref{example:02} under some strong assumptions, which provide intuitive insight regarding the much more general theory presented in Sections~\ref{sec:02} and~\ref{multiple-sike-models}.

\begin{figure}[h]
 \begin{center}
 \includegraphics[width=\textwidth]{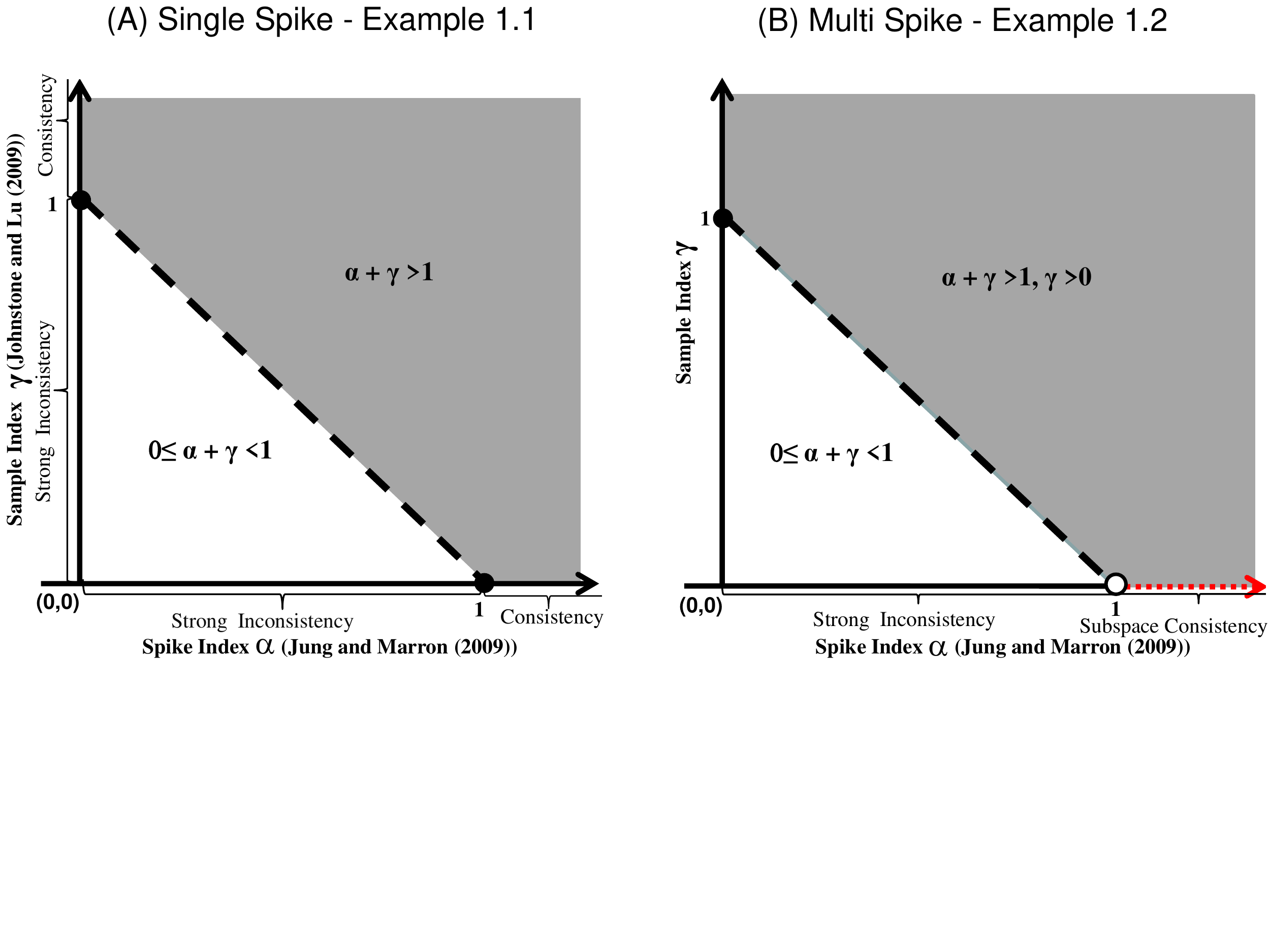}
 \end{center}
 \vspace{-3cm}
 \caption{General consistency and strong inconsistency regions for PCA, as a function of the spike index $\alpha$ and
 the sample index $\gamma$. 
{{Panel (A) - single spike model in Example~\ref{example:01}:  PCA is consistent on
 the grey region ($\alpha+\gamma>1$),
 and strongly inconsistent on the white triangle ($0 \leq \alpha+\gamma<1$).
Panel (B) - multiple spike model in Example~\ref{example:02}: the first $m$ sample PCs 
are consistent on the grey region ($\alpha+\gamma>1, \gamma>0$), subspace consistent on the dotted line segment ($\alpha>1, \gamma=0$) on the horizontal}} axis,
and strongly inconsistent 
on the white triangle ($0 \leq \alpha+\gamma<1$). 
 }\label{fig:01}
 \vspace{-.3cm}
 \end{figure}

For these two illustrative examples, the three types of information and their relationships can be mathematically quantified by two indices, namely the {\it spike index} $\alpha$ and the {\it sample index} $\gamma$. Within the context of these examples, we point out the significant contributions of our results in comparison with existing results. The comparisons and connections are graphically illustrated in Figure~\ref{fig:01} and discussed below. 

\begin{example} (Single-component spike model) Assume that $X_1,\ldots, X_n$ are random sample vectors
from a $d$-dimensional normal distribution $N(0, \Sigma)$, where the sample size
$n \sim d^\gamma$ ($\gamma\geq 0$ is defined as the sample index)
and the covariance matrix $\Sigma$
has the eigenvalues as
\begin{equation*}
\lambda_1 \sim d^\alpha, \lambda_2=\cdots=\lambda_d=1, \alpha\geq0,
\end{equation*}
where the constant $\alpha$
is defined as the spike index.

Theorem~\ref{Th:01}, when applied to this example, suggests that the maximal sample eigenvector
is consistent when $\alpha+\gamma>1$ (grey region in Figure~\ref{fig:01}(A)), and strongly inconsistent when
$0 \leq\alpha+\gamma<1$ (white triangle in Figure~\ref{fig:01}(A)).
 These very general new results nicely connect with many
existing ones:
\begin{itemize}

\item \textbf{Previous Results I - the classical domain}:

For this example, Theorem 1 of Anderson~\cite{anderson1963asymptotic} implied that for fixed dimension $d$ and finite eigenvalues,
 when the sample size $n \rightarrow \infty$ (i.e. $\gamma \rightarrow \infty$, the limit on the vertical axis), the maximal sample eigenvector is consistent. This case is the upper left corner of Figure~\ref{fig:01}(A).

\item \textbf{Previous Results II - the random matrix domain}:\begin{enumerate}
\item[(a)] The results of Johnstone and Lu~\cite{johnstone2009consistency} appear
 on the vertical axis in Panel (A)  where the spike index $\alpha=0$ (as they fix the spike information):
the first sample eigenvector
is consistent
when the {sample index} $\gamma>1$ and strongly inconsistent when $\gamma<1$.
\item[(b)] {{Nadler~\cite{nadler2008finite} explored the interesting boundary case of $\alpha=0, \gamma=1$ (i.e. $\frac{d}{n}\rightarrow c$ for a constant $c$) and
    showed that $<\hat{u}_1, u_1>^2 \xrightarrow{\rm a.s}\frac{((\lambda_1-1)^2-c)_{+}}{(\lambda_1-1)^2+c(\lambda_1-1)}$,
    where $\hat{u}_1$ and $u_1$ are the first sample and population eigenvector.}}
     This result appears in Panel (A) as the single solid circle $\gamma=1$ on the vertical axis.
\end{enumerate}

%

\item \textbf{Previous Results III - the HDLSS domain}:
\begin{enumerate}
\item[(a)]
 The theorems of Jung and Marron~\cite{jung2009pca} are represented on the horizontal axis in Panel (A) when the sample index $\gamma=0$ (as they fix the sample size):
the maximal sample eigenvector
 is consistent with
the first population eigenvector when the {spike index} $\alpha>1$ and strongly
inconsistent when $\alpha<1$.

\item[(b)] {{Jung et al.~\cite{Jung2010} deeply explored limiting behavior at the boundary $\alpha=1, \gamma=0$
(i.e. $\frac{d}{\lambda_1}\rightarrow c$ for a constant $c$) and showed that
$<\hat{u}_1, u_1>^2 \Rightarrow\frac{\chi_n^2}{\chi_n^2+ c}$,
 where $``\Rightarrow"$ means convergence in distribution and $\chi_n^2$ is
 the chi-squared distribution with $n$ degrees of freedom.
This result appears in Panel (A) as the single solid circle $\alpha=1$ on the horizontal axis.}}

\end{enumerate}
\item \textbf{Our Results} {{hence nicely connect existing domains of asymptotics,
and give a much more complete characterization for the regions of PCA consistency, subspace consistency, and strong inconsistency.   We also investigate asymptotic properties of the other sample eigenvectors and  all the sample eigenvalues.}}
\end{itemize}
\label{example:01}
\end{example}

\begin{example} (Multiple-component spike model) Assume that the covariance matrix $\Sigma$ in Example~\ref{example:01}
has the following eigenvalues
\begin{equation*} \lambda_j=
\begin{cases} c_j d^\alpha &\text{if $j\leq m $,}
\\
1 &\text{if $j> m $,}
\end{cases}
\quad {\alpha\geq0},
\end{equation*}
where $m$ is a finite positive integer, the constants $c_j, j=1,\cdots, m$, are positive and satisfy that $c_j>c_{j+1}>1$, $j=1,\cdots, m-1$.

Theorem~\ref{Th:02}, when applied to this example, shows that the first $m$ sample eigenvectors
 are individually consistent with corresponding population eigenvectors when $\alpha+\gamma>1, \gamma>0$ (the grey region in Figure~\ref{fig:01}(B)), instead of being subspace consistent~\citep{jung2009pca}, and strongly inconsistent when
$\alpha+\gamma<1$ (the white triangle in Panel (B)).
This very general new result connects with many
others in the existing literature:

\begin{itemize}

\item \textbf{Previous Results I - the classical domain}:

For this example, Theorem 1 of Anderson~\cite{anderson1963asymptotic} implied that for fixed dimension $d$ and finite eigenvalues, 
when the sample size $n \rightarrow \infty$ (i.e. $\gamma \rightarrow \infty$, the limit on the vertical axis), the first $m$ sample eigenvectors
  are consistent,  
while the other sample eigenvectors are
 subspace consistent. 
 This case is the upper left corner of Figure~\ref{fig:01}(B).

\item \textbf{Previous Results II - the random matrix domain}:

{{Paul~\cite{paul2007asymptotics} explored asymptotic properties of the first
$m$ eigenvectors and eigenvalues in the interesting boundary case of $\alpha=0, \gamma=1$, i.e., $\frac{d}{n}\rightarrow c$ with $c\in(0,1)$ and showed that $<\hat{u}_j, u_j>^2 \xrightarrow{\rm a.s}\frac{((\lambda_j-1)^2-c)_{+}}{(\lambda_j-1)^2+c(\lambda_j-1)}$ for
$j=1,\cdots, m$.
This result appears in Panel (B) as the solid circle $\gamma=1$ on the vertical axis.}}
Paul and Johnstone~\cite{paul2007Davis} 
considered a similar
 framework but from a minimax risk analysis perspective. Nadler~\cite{nadler2008finite} and Johnstone and Lu~\cite{johnstone2009consistency} did not study multiple spike models.

\item \textbf{Previous Results III - the HDLSS domain}:

The theorems of Jung and Marron~\cite{jung2009pca} are valid on  the horizontal axis in Panel (B) where the sample index $\gamma=0$. In particular, for this example, their results
showed
 that the first $m$ sample eigenvectors
 are not separable when the spike index $\alpha>1$ (the horizontal dotted red line segment), instead they are {subspace consistent} with
their corresponding  population eigenvectors,
 and are {strongly inconsistent} when the spike index $\alpha<1$ (the horizontal solid line segment).
 They and Jung et al.~\cite{Jung2010}
 did not study the asymptotic behavior on the boundary - the single open circle $(\alpha=1, \gamma=0)$ on the horizontal axis.

\item \textbf{Our Results}
cover the classical domain,
 and are stronger than what~\cite{jung2009pca} obtained:
 the increasing sample size enables us to separate out the first few leading eigenvectors and characterize {\it individual consistency},
 while only subspace consistency was obtained by~\cite{jung2009pca}. 
\end{itemize}
\label{example:02}
\end{example}

The organization of the rest of the paper is as follows. Section~\ref{subsec:11} first introduces our notations and several relevant consistency concepts. 
Section~\ref{sec:02} then presents the
theoretical results of single-component spike models, stating the asymptotic properties
of the sample eigenvalues and eigenvectors under our general framework.
Section~\ref{ssec:single:increasingn} first considers single-component spike models with the increasing
sample size $n$,
and Section~\ref{ssec:single:fixedn} then studies
single-component spike models where the sample size $n$ is fixed.
Section~\ref{multiple-sike-models} studies multiple-component spike models. For easy access to the main ideas, Section~\ref{ssec:22} first studies models  with distinct eigenvalues, while Section~\ref{ssec:23}  then considers models where the eigenvalues are grouped. {{Section~\ref{dicussion} contains some discussion about the asymptotic properties of PCA when
some small eigenvalues equal to zero and
the challenges to obtain non-asymptotic results.}}
Section~\ref{proofs} contains the technical proofs of the main theorem.

\section{Notations and Concepts}\label{subsec:11}
We now introduce some necessary notations, and define consistency concepts relevant for our asymptotic study. 

\subsection{Notation}\label{subsec:notation}
Let the population covariance matrix be $\Sigma$, whose eigen decomposition is
$$
\Sigma=U\Lambda U^T,
$$ where $\Lambda$ is the diagonal matrix of population eigenvalues
$\lambda_1\geq\lambda_2\geq\ldots\geq\lambda_d$, and $U$
is the  matrix of corresponding eigenvectors
$U=[u_1,\ldots,u_d]$.


{{As in Jung and Marron~\cite{jung2009pca},  assume that $X_1,\ldots, X_n$ are i.i.d.
$\;d$-dimensional random sample vectors and have the following representation
\begin{equation}~\label{subguassion:assumption}
X_i=\sum_{j=1}^d z_{i,j} u_j,
\end{equation}
where the $z_{i,j}$'s are i.i.d random variables with zero mean, unit variance and
finite fourth moment. An important special case is that the $z_{i,j}$'s follow
 the standard normal distribution $N(0,1)$.


{{\begin{assumption}
$X_1,\ldots, X_n$ are a
random sample having the distribution described by ~\eqref{subguassion:assumption}.
\label{gauss-assumption}
\end{assumption}
Denote the sample covariance matrix by $\hat{\Sigma}=n^{-1}XX^{T}$, where $X=[X_1,\ldots,X_n]$.
Note that $\hat{\Sigma}$ can also be
decomposed as}}
\begin{equation}
\hat{\Sigma}=\hat{U}\hat{\Lambda}\hat{U}^T,
\label{sample:covariance}
\end{equation}
where $\hat{\Lambda}$ is the diagonal matrix of sample eigenvalues
$\hat{\lambda}_1\geq\hat{\lambda}_2\geq\ldots\geq\hat{\lambda}_d$
and $\hat{U}$ is the matrix of corresponding sample eigenvectors where
$\hat{U}=[\hat{u}_1,\ldots,\hat{u}_d]$.

{{Below we introduce asymptotic notations that will be used in our theoretical studies.
Assume that $\{\xi_{n}: n=1,\ldots,\infty\}$
 is a sequence of random variables,  and
$\{a_{n}: n=1,\ldots,\infty\}$ is a sequence of constant values.
\begin{itemize}
\item Denote $\xi_{n}={\rm o}_{\rm a.s}\left(a_{n}\right)$ if $\mbox{lim}_{n\rightarrow \infty}\frac{\xi_{n}}{a_{n}}=0$ almost surely.
\item Denote $\xi_{n}={\rm O}_{\rm a.s}\left(a_{n}\right)$ if $\overline{\mbox{lim}}_{n\rightarrow \infty}\left|\frac{\xi_{n}}{a_{n}}\right|\leq z$ almost surely,
where the random variable $z$ satisfies $P(0 < z < \infty)=1$..
\item Denote $\xi_{n}\stackrel{\rm a.s}{\sim}a_{n}$ if $c_2 \leq \underline{\mbox{lim}}_{n\rightarrow \infty}\frac{\xi_{n}}{a_{n}}\leq
\overline{\mbox{lim}}_{n\rightarrow \infty}\frac{\xi_{n}}{a_{n}}\leq c_1$ almost surely,
for two constants $c_1\geq c_2>0$.
\end{itemize}}}


In addition, we introduce the following notions to help understand the assumptions on the population
 eigenvalues in our theorems.
Assume that $\{a_k: k=1,\ldots,\infty\}$ and $\{b_k: k=1,\ldots,\infty\}$ are two sequence of constant values, where
$k$ can stand for either $n$ or $d$.
\begin{itemize}

\item Denote $a_{k} \gg b_{k}$ if $\mbox{lim}_{k\rightarrow \infty}\frac{b_{k}}{a_{k}}=0$.
\item Denote $a_{k}\sim b_{k}$ if $c_2 \leq
\underline{\mbox{lim}}_{k\rightarrow \infty}\frac{a_k}{b_k}\leq
\overline{\mbox{lim}}_{k\rightarrow \infty}\frac{a_k}{b_k}\leq c_1$ for two constants $c_1\geq c_2>0$.
\end{itemize}

\subsection{Concepts}\label{subsec:concept}
Below we list three important concepts relevant for consistency and strong inconsistency, some of
which are modified from the related concepts given by
Jung and Marron~\cite{jung2009pca} and Shen et al.~\cite{shen2011}.


Let $\hat{u}_j$ be
any normalized sample estimator of $u_j$ for $j=1,\ldots ,
 [n\wedge d]$. 


 \begin{itemize}

\item \textbf{Consistency with rate $a_n$}:
 The estimator $\hat{u}_j$ is consistent with
its population counterpart $u_j$ with the convergence rate $a_n$ if
$\mid<\hat{u}_j,u_j>\mid=1+{\rm O}_{\rm a.s}(a_n)$. For example,
$a_n=\left(\frac{n\lambda_1}{d}\right)^{\frac{1}{2}}$.

\item \textbf{Strong inconsistency with rate $a_n$}:
$\hat{u}_j$ is strongly inconsistent with
 $u_j$ with the convergence rate $a_n$ if
$|<\hat{u}_j,u_j>|={\rm O}_{\rm a.s}(a_n)$.
\end{itemize}

Let $H$ be an index set, e.g. $H=\{m+1, \cdots, d\}$. Define $S=\mbox{span}\{u_k, k\in H\}$ to
 be the linear span generated
by $\{u_k, k\in H\}$.
\begin{itemize}
\item \textbf{Subspace consistency with rate $a_n$}:
$\hat{u}_j$, $j\in H$, is subspace consistent with
$S$ with convergence rate $a_n$ if
\begin{eqnarray}
\mbox{angle}(\hat{u}_j, S)
={\rm O}_{a.s}(a_n),
\label{distance}
\end{eqnarray}
where the angle between the estimator $\hat{u}_j$ and the subspace $S$ is the
angle between the estimator and its projection onto the subspace, see Jung and Marron~\cite{jung2009pca}.
For further clarification, we provide a graphical illustration of the angle in Section B of the supplement~\cite{shen2011online}.

\end{itemize}

{{Terminology: In the following results, for simple general formulations,
the term - \emph{consistent at the rate $a_n$} - will mean $``{\rm O}_{\rm a.s}(a_n)"$
in situations where $a_n\rightarrow 0$.  Otherwise it means $``{\rm o}_{\rm a.s}(a_n)"$.
 Similarly for \emph{strong inconsistency} and \emph{subspace consistency}.}}



\section{Single component spike models}\label{sec:02}
Below we state our main theorems for single-component spike models.
In Section~\ref{ssec:single:increasingn}, we study the asymptotic properties of PCA
with increasing sample size $n$. 
In Section~\ref{ssec:single:fixedn} , we investigate the asymptotic properties of PCA with fixed $n$.



\subsection{Cases with increasing sample size $n$}~\label{ssec:single:increasingn}
We first state in Theorem~\ref{Th:01} one of our main theoretical results regarding PCA consistency under our general framework. We then offer several remarks in regards to the conditions of the theorem as well as the connection between our results and the earlier ones in the literature.

To fix ideas, we assume the maximal eigenvalue $\lambda_1$ dominates the other eigenvalues.  WLOG, we assume that as $n\rightarrow \infty$ or $d\rightarrow \infty$,
\begin{assumption}
$\lambda_1> \lambda_2\rightarrow \cdots\rightarrow \lambda_d\rightarrow c_{\lambda}$, where $c_\lambda$ is a constant.
\label{eigenvalue-assumption-single}
\end{assumption}
 As discussed in the Introduction, we consider the delicate balance among the positive {\it sample information} $n$, the positive {\it spike information} $\lambda_1$, and the negative {\it variable information} $d$, and characterize the various PCA
 consistency and strong-inconsistency regions.

{{Theorem~\ref{Th:01} below suggests that the asymptotic properties of the sample eigenvalues and eigenvectors depend on the relative strength of the positive information and the negative information, as particularly measured by two ratios: $\frac{d}{n\lambda_1}$ and $\frac{d}{n}$. The value of $\frac{d}{n\lambda_1}$ determines whether the maximal sample eigenvalue is separable from the other eigenvalues, and further determines the consistency of the maximal sample  eigenvector.
The value of $\frac{d}{n}$ determines the asymptotic properties of the second and higher sample eigenvalues and
 eigenvectors.}}

The following discussion and the scenarios in Theorem~\ref{Th:01} are arranged
according to a decreasing amount of positive information:
{{\begin{itemize}
\item Theorem~\ref{Th:01}(a): If the amount of positive information dominates the amount of negative information up to the maximal eigenvalue, i.e. $\frac{d}{n\lambda_1}\rightarrow 0$, then the maximal sample eigenvector is consistent, and the other sample eigenvectors are subspace consistent.
    In addition, the asymptotic properties of sample eigenvalues and eigenvectors whose index are greater than 1 depend on the value $\frac{d}{n}$.
\item Theorem~\ref{Th:01}(b): On the other hand, if the amount of negative information always dominates, i.e. $\frac{d}{n\lambda_1}\rightarrow \infty$, then the sample eigenvalues are asymptotically indistinguishable, and the sample eigenvectors are strongly inconsistent.
\end{itemize}}}




\begin{theorem}  Under Assumptions~\ref{gauss-assumption} and~\ref{eigenvalue-assumption-single},
as $n\rightarrow \infty$, the following results hold.

\begin{enumerate}

\item[(a)] If $\frac{d}{n\lambda_1}\rightarrow 0$, then $\frac{\hat{\lambda}_1}{\lambda_1}\xrightarrow{\rm a.s} 1$,
$\hat{u}_1$ is consistent with $u_1$,
 and the other $\hat{u}_j$ are subspace consistent
with $S=\mbox{span}\{u_k, k\geq2\}$. 
In addition,
\begin{enumerate}
\item[i.] {{If $\frac{d}{n}\rightarrow 0$, then $\frac{\hat{\lambda}_j}{\lambda_j}\xrightarrow{\rm a.s} 1$, $j=2,\cdots, [n\wedge (d-1)]$,
  and $\hat{\lambda}_{[n\wedge d]}={\rm O}_{\rm a.s}(1)$. The consistency rate for $\hat{u}_1$ and
  the subspace consistency
  rate for the other $\hat{u}_j$ are both $\left(\frac{1}{\lambda_1}\right)^{\frac{1}{2}}$.}}

\item[ii.] {{If $\frac{d}{n}\rightarrow \infty$, then $\frac{\hat{\lambda}_j}{\lambda_j} \stackrel{\rm a.s}{\rightarrow} \frac{d}{n} $
for $j=2,\cdots, [n\wedge d]$;
$\hat{u}_1$ is consistent with rate
$\left(\frac{d}{n\lambda_1}\right)^{\frac{1}{2}}$, and the other $\hat{u}_j$ are strongly
inconsistent with rate $\left(\frac{n\lambda_j}{d}\right)^{\frac{1}{2}}$.}}

\item[iii.] {{If  $\overline{{\rm lim}}\frac{d}{n}=c^*$ ($0<c^*\leq \infty$), then
 for $j=2,\cdots, [n\wedge d]$,
 $\overline{{\rm lim}}\hat{\lambda}_j \leq c \times \overline{{\rm lim}}\frac{d}{n}$ almost surely,
 where $c$ is some constant. The consistency rate for $\hat{u}_1$ and the subspace consistency
  rate for the other $\hat{u}_j$ are both $\left(\frac{d_n}{n\lambda_1}\right)^{\frac{1}{2}}$,
 where $\{d_n\}$ is a sequence converging to $c^*$.}}
\end{enumerate}


\item[(b)] If $\frac{d}{n\lambda_1}\rightarrow \infty$,
then  $\hat{\lambda}_j \stackrel{\rm a.s}{\rightarrow}c_\lambda \frac{d}{n}$
for $j=1,\cdots, [n\wedge d]$,
and the corresponding eigenvectors $\hat{u}_j$
are strongly inconsistent with rate $(\frac{n\lambda_j}{d})^{\frac{1}{2}}$.

\end{enumerate}
\label{Th:01}
\end{theorem}

Having stated the main results for single-component spike models, we now offer several remarks regarding
the conditions assumed in Theorem~\ref{Th:01}
and make connections with  existing results about PCA consistency.

\begin{itemize}

%
%
%
%
%
%
%
%
%
%
%
%
%
%

\item {{If Assumption~\ref{eigenvalue-assumption-single} is replaced by the alternative assumption
$\lambda_1\gg \lambda_2\sim \cdots\lambda_d\sim 1$, then except for  $\frac{\hat{\lambda}_1}{\lambda_1}\xrightarrow{\rm a.s} 1$
in Scenario (a), all other $``\xrightarrow{\rm a.s}"$ for the sample eigenvalues should be replaced by $``\stackrel{\rm a.s}{\sim}"$.
The results for the sample eigenvectors remain the same.}}

\item An assumption of the form (3.1), i.e $\lambda_2\rightarrow \cdots\rightarrow \lambda_d\rightarrow c_{\lambda}$, or
else $\lambda_2\sim \cdots\lambda_d\sim 1$ is needed to obtain general convergence results for the non-spike sample eigenvalues
$\hat{\lambda}_j$, $j>1$ under the  wide range of  scenarios: $\frac{d}{n}\rightarrow 0$,
 $\frac{d}{n}\rightarrow \infty$ or $\overline{{\rm lim}}\frac{d}{n}=c^*$ ($0<c^*\leq \infty$). {When one focusses only on the spike eigenvalue, a weaker assumption, such as the slowly decaying non-spike eigenvalues assumed by
 Bai and Yao (2012)~\cite{bai2012sample}, is enough.
Then the spike condition $\lambda_1\gg \lambda_2$ is enough to generate the consistency properties of $\hat{\lambda}_1$ and
 $\hat{u}_1$ in Scenario (a).
In that case, the behaviors of the other sample
 eigenvalues and  eigenvectors are very case-wise to formulate in general.}

\item  {{Nadler~\cite{nadler2008finite} and
Johnstone and Lu (2009)~\cite{johnstone2009consistency} studied the properties of the first
sample eigenvalue and eigenvector under the normality assumption for  fixed $\lambda_1$.
Furthermore, if $d$ is fixed, {Scenario (a) of Theorem~\ref{Th:01}} degenerates to the case
 studied by Anderson (1963)~\cite{anderson1963asymptotic}.}}

\item Assuming fixed $\lambda_1$ and $\frac{d}{n}\rightarrow c$ with $c$ being a constant, Nadler~\cite{nadler2008finite}, Johnstone and Lu~\cite{johnstone2009consistency}
and Benaych-Georges and Nadakuditi~\cite{benaych2011eigenvalues} obtained the results in
{\it Previous Results II - the random matrix domain} in Example~\ref{example:01},
which indicate that, as
$n\rightarrow \infty$, the maximal sample eigenvector
 $\hat{u}_1$ is consistent when $\frac{d}{n}\rightarrow 0$, and
inconsistent when $\frac{d}{n}\rightarrow \infty$. 
Our Theorem~\ref{Th:01} includes this as a special case.
~In addition, Theorem~\ref{Th:01} offers more than just relaxing the fixed $\lambda_1$ assumption: it characterizes how an increasing $\lambda_1$ interacts with the ratio $\frac{d}{n}$, derives the corresponding convergence rate, and also studies the asymptotic properties of the higher order sample eigenvalues and eigenvectors, all of which have not been investigated before.

\end{itemize}

\subsection{Cases with fixed $n$}~\label{ssec:single:fixedn}
 Theorem~\ref{Th:01:hdlss} summarize the results for the fixed $n$ cases (i.e. the HDLSS domain). In comparison with Jung and Marron (2009)~\cite{jung2009pca}, we make more general assumptions on the population eigenvalues, and obtain the corresponding convergence rate results; furthermore, we obtain almost sure convergence, instead of convergence in probability~\cite{jung2009pca}.

Consider the $z_{i,j}$ in~\eqref{subguassion:assumption}, and define
\begin{equation}
 \widetilde{Z}_j=(z_{1,j},\cdots, z_{n,j})^T, \quad j=1,\cdots, d,
\label{dual:Z}
 \end{equation}
 which are needed here to describe the asymptotic
  properties of the sample eigenvalues in HDLSS settings. In addition, define $K=\lim_{d\rightarrow \infty}\frac{\sum_{j=2}^d \lambda_j}{nd}$.

\begin{theorem}
 Under Assumptions~\ref{gauss-assumption} and~\ref{eigenvalue-assumption-single},
 for fixed $n$, as $d\rightarrow \infty$, the following results hold.

\begin{enumerate}
\item[(a)] If $\frac{d}{\lambda_1}\rightarrow 0$,
then  $\frac{\hat{\lambda}_1}{\lambda_1 }\stackrel{\rm a.s}{\rightarrow}
\frac{\widetilde{Z}_1^T \widetilde{Z}_1}{n}$,
where $\widetilde{Z}_1$ is defined in~\eqref{dual:Z},
and the rest of the
non-zero
$\frac{\hat{\lambda}_j}{d} \stackrel{\rm a.s}{\rightarrow} K$.
In addition,
$\hat{u}_1$ is consistent  with rate
$\left(\frac{d}{\lambda_1}\right)^{\frac{1}{2}}$, and the rest of the $\hat{u}_j$
are strongly inconsistent with rate
$\left(\frac{\lambda_j}{d}\right)^{\frac{1}{2}}$.

\item[(b)] If $\frac{d}{\lambda_1}\rightarrow \infty$,
then the non-zero $\frac{\hat{\lambda}_j}{d}\stackrel{\rm a.s}{\rightarrow} K$, and the corresponding
 $\hat{u}_j$ are strongly inconsistent
 with rate $(\frac{\lambda_j}{d})^{\frac{1}{2}}$,
 respectively.

\end{enumerate}
\label{Th:01:hdlss}
\end{theorem}

Some comments about the conditions  and results of Theorem~\ref{Th:01:hdlss}
\begin{itemize}
\item {{Assumption~\ref{eigenvalue-assumption-single} can be replaced by
$\lambda_1\gg \lambda_2\sim \cdots\lambda_d\sim 1$. The results remain the same.}}

\item {{Even if the non-spike eigenvalues $\lambda_j$, $j>2$, decay slowly,
 the condition $\lambda_1\gg \lambda_2$ is enough to generate the same properties for $\hat{\lambda}_1$ and
 $\hat{u}_1$ as in Scenario (a).}}

\item {{If Assumption~\ref{gauss-assumption} is strengthened to a normality assumption, then
$\frac{\hat{\lambda}_1}{\lambda_1}  \stackrel{\rm a.s}{\rightarrow} \frac{\chi_n^2}{n}$ in Scenario (a).}}

\item {{Assumption~\ref{gauss-assumption} assumes that the $z_{i,j}$'s are i.i.d rather than $\rho$-mixing as
in~\cite{jung2009pca}.
Thus, convergence in probability in~\cite{jung2009pca} is strengthened to almost sure convergence here.}}

\end{itemize}

\section{Multiple component spike models}
\label{multiple-sike-models}
We consider multiple spike models with {{finite}} $m(\in [1,n\wedge d])$ dominating spikes. In Section~\ref{ssec:22}, we study models where the dominating eigenvalues are distinct.
 In Section~\ref{ssec:23},  we consider the cases where the eigenvalues are not all distinct,
 by introducing the concept of \emph{tiered eigenvalues}.

\subsection{Multiple component spike models with distinct eigenvalues}
\label{ssec:22}

\subsubsection{Cases with increasing sample size $n$} {{WLOG, we assume that
 the first $m$ population eigenvalues have different strength and dominate the rest population eigenvalues,
 which are asymptotically equivalent. }}
\begin{assumption}~\label{eigenvalue-assumption-multiple}
as $n \rightarrow \infty$, $\lambda_1 > \cdots > \lambda_m> \lambda_{m+1}\rightarrow \cdots\rightarrow \lambda_d\rightarrow
c_\lambda.$
\end{assumption}
A useful quantity, for distinguishing the various cases among eigenvectors in the coming theorems, is
\begin{equation*}
a_l=\mbox{max}_{1\leq k \leq l}\frac{\lambda_{k+1}}{\lambda_k}, \quad l=1,\cdots,m.
\end{equation*}
{{This lower bound on the consecutive relative gap among the first $l$ eigenvalues provides a critical measure of the separation between the $l$-th sample eigenvector and the first $l-1$ sample eigenvectors.}}

Below we first state the main theoretical results in Theorem~\ref{Th:02}, and follow up with some remarks about the theorem conditions and the connections between the theorem and the existing results in the literature.

Similar to Theorem~\ref{Th:01}, Theorem~\ref{Th:02} states the asymptotic properties of the sample eigenvalues and eigenvectors in a trichotomous manner, separated by the size of $\frac{d}{n\lambda_j}$, which again measures the relative strength of the positive information and the negative information. The three scenarios below and in Theorem~\ref{Th:02}  are arranged in a decreasing order of the amount of the positive information:
\begin{itemize}
\item Theorem~\ref{Th:02}(a): {{If the amount of positive information dominates the amount of negative information up to the $m$th spike, i.e. $\frac{d}{n\lambda_m}\rightarrow 0$,
then each of the first $m$ sample eigenvector is consistent, and the additional ones are subspace consistent;}}
\item Theorem~\ref{Th:02}(b): {{Otherwise, if the amount of positive information dominates the amount of negative information only up to the $h$th spike $(h\in[1,m])$, i.e. $\frac{d}{n\lambda_h}\rightarrow 0$ and $\frac{d}{n\lambda_{h+1}}\rightarrow \infty$, then each of the first $h$ sample eigenvector is consistent, and each of the remaining higher-order sample eigenvector is strongly-inconsistent;}}
\item Theorem~\ref{Th:02}(c): Finally, if the amount of negative information always dominates, i.e. $\frac{d}{n\lambda_1}\rightarrow \infty$, then the sample eigenvalues are asymptotically indistinguishable, and the sample eigenvectors are strongly inconsistent.
\end{itemize}


\begin{theorem}
  Under Assumptions~\ref{gauss-assumption} and~\ref{eigenvalue-assumption-multiple}, as $n \rightarrow \infty$, the following results hold.

\begin{enumerate}

\item[(a)] \emph{If $\frac{d}{n\lambda_m}\rightarrow 0$, then $\frac{\hat{\lambda}_j}{\lambda_j} \xrightarrow{\rm a.s} 1$
 for $1\leq j\leq m$.
 In addition, $\hat{u}_j$ are consistent with $u_j$
for $1\leq j\leq m$ and the other $\hat{u}_j$ are subspace consistent with $S=\mbox{span}\{u_k, k\geq m+1\}$.}

\item[(b)] If there exists a constant $h$, $1\leq h \leq m$, such that
$\frac{d}{n\lambda_h}\rightarrow 0$ and $\frac{d}{n\lambda_{h+1}}\rightarrow \infty$,
 then $\frac{\hat{\lambda}_j}{\lambda_j} \xrightarrow{\rm a.s} 1$ for $1\leq j\leq h$,
 and the other non-zero $\hat{\lambda}_j \stackrel{\rm a.s}{\rightarrow} c_\lambda \frac{d}{n} $.
 In addition, $\hat{u}_j$ are consistent with rate
$\left(a_j\vee \frac{d}{n\lambda_j}\right)^{\frac{1}{2}}$  for $1\leq j\leq h$,
and the other $\hat{u}_j$ are strongly inconsistent with
rate $\left(\frac{n\lambda_j}{d}\right)^{\frac{1}{2}}$.

\item[(c)] If $\frac{d}{n\lambda_1}\rightarrow \infty$,
then the non-zero $\hat{\lambda}_j \stackrel{\rm a.s}{\rightarrow}c_\lambda\frac{d}{n} $, and the corresponding
 $\hat{u}_j$ are strongly inconsistent with rate $(\frac{n\lambda_j}{d})^{\frac{1}{2}}$.

\end{enumerate}
\label{Th:02}
\end{theorem}

We now discuss the properties of the rest of the sample eigenvalues and
the convergence rate in Scenario (a), and the conditions needed in the theorem
 and how the results connect with  existing ones in the literature.

\begin{itemize}

\item {The special case of $m=1$ is  Theorem~\ref{Th:01} for single spike models.}

\item As in Theorem~\ref{Th:01},  Scenario (a) in Theorem~\ref{Th:02} contains three different cases.

\begin{enumerate}

 \item[i.] If $\frac{d}{n}\rightarrow 0$, then
 $\frac{\hat{\lambda}_j}{\lambda_j}\xrightarrow{\rm a.s} 1$, $j=m+1,\cdots, [n\wedge (d-m)]$
 and the rest of the non-zero $\hat{\lambda}_j={\rm O}_{\rm a.s}(1)$. In addition, the consistency rates for
 the $\hat{u}_j$ are  $\left(a_j\vee \frac{1}{\lambda_j}\right)^{\frac{1}{2}}$
  for $1\leq j\leq m$, and   the subspace consistency rates for the other $\hat{u}_j$ are
  $\left(a_m\vee \frac{1}{\lambda_m}\right)^{\frac{1}{2}}$.


 \item[ii.] {The case $\frac{d}{n}\rightarrow \infty$  is considered in
  Scenario (b) ($h=m$) of Theorem~\ref{Th:02}.}

 \item[iii.] If $\overline{{\rm lim}}\frac{d}{n}=c^*$ ($0<c^*\leq \infty$), then
for $j=m+1,\cdots, [n\wedge d]$,
$\overline{{\rm lim}}\hat{\lambda}_j \leq c \times \overline{{\rm lim}}\frac{d}{n}$ almost surely,
 where $c$ is some constant. 
 Also the consistency rates for the $\hat{u}_j$ are  $\left(a_j\vee \frac{d_n}{n\lambda_j}\right)^{\frac{1}{2}}$
  for $1\leq j\leq m$, and the  subspace consistency rates for the rest of the $\hat{u}_j$ are
  $\left(\frac{d_n}{n\lambda_m}\right)^{\frac{1}{2}}$, where $\{d_n\}$ is a sequence converging to
   $c^*$.


\end{enumerate}

\item {If Assumption~\ref{eigenvalue-assumption-multiple} is replaced by the alternative assumption
$\lambda_1>\cdots >\lambda_m \gg \lambda_{m+1}\sim \cdots\lambda_d\sim 1$,
then we still have $\frac{\hat{\lambda}_j}{\lambda_j}\xrightarrow{\rm a.s} 1$, $1\leq j\leq m$, as in Scenario (a) and
$1\leq j\leq h$ as in Scenario  (b),
but all other results of the form $``\xrightarrow{\rm a.s}"$  for the sample eigenvalues should be replaced by $``\stackrel{\rm a.s}{\sim}"$.
The results for the sample eigenvectors remain same.}

\item {Even if the non-spike eigenvalues $\lambda_j$, $j>m$, decay slowly,
 the condition $\lambda_1>\cdots >\lambda_m \gg \lambda_{m+1}$ is enough to generate the consistency properties of $\hat{\lambda}_j$ and
 $\hat{u}_j$, for $1\leq j\leq m$ in Scenario (a) and $1\leq j\leq h$ in Scenario (b).}


\item {In Theorem~\ref{Th:02}, consider the special case of fixed dimension $d$ and
$\infty >\lambda_1 > \cdots > \lambda_m>  \lambda_{m+1}\rightarrow\cdots\rightarrow \lambda_d \rightarrow c_\lambda$.
Then, Theorem~\ref{Th:02}(a) is consistent with
 the classical results implied by Theorem 1 of Anderson~\cite{anderson1963asymptotic}.}



\item Considering fixed $\lambda_1,\cdots, \lambda_m$ and $\frac{d}{n}\rightarrow c$, where  $c \in (0,1)$, Paul~\cite{paul2007asymptotics}
 obtained results that are applicable to Example~\ref{example:02} to obtain {\it Previous Results II - the random matrix domain} in .
  As one can see, our Theorem~\ref{Th:02} relaxes the assumptions
   of $\frac{d}{n}\rightarrow c \in (0,1)$ and that $\lambda_1, \cdots, \lambda_m$ are fixed. In addition,
    we characterize how increasing $\lambda_1,\cdots, \lambda_m$ interact with the ratio
    $\frac{d}{n}$ along with the corresponding convergence rates,
    and study the asymptotic properties of the higher order sample eigenvalues and eigenvectors,
    all of which have not been investigated before.

\end{itemize}

\subsubsection{Cases with fixed $n$} The following Theorem~\ref{Th:02:hdlss} considers cases with fixed $n$. {{The multiple spike condition in Assumption~\ref{eigenvalue-assumption-multiple} now becomes
that the first $m$ population eigenvalues are of the different order and dominate the other population eigenvalues,
 which are asymptotically equivalent:}}
\begin{assumption}
as $d \rightarrow \infty$, $\lambda_1 \gg \cdots \gg \lambda_m\gg \lambda_{m+1}\sim\cdots\sim \lambda_d\sim 1.$
\label{eigenvalue-assumption-multiple-hdlss}
\end{assumption}
Note that for fixed $n$ and $d\rightarrow \infty$, assuming $\lambda_j>\lambda_{j+1}$ can not
asymptotically separate the corresponding sample eigenvalues $\hat{\lambda}_j$ and $\hat{\lambda}_{j+1}$.
Thus, we need to replace Assumption~\ref{eigenvalue-assumption-multiple}
with Assumption~\ref{eigenvalue-assumption-multiple-hdlss} to asymptotically separate
the first $m$ sample eigenvalues.
Define $K=\lim_{d\rightarrow \infty} \frac{\sum_{j=m+1}^d \lambda_j}{nd}$.

\begin{theorem}  Under Assumptions~\ref{gauss-assumption} and~\ref{eigenvalue-assumption-multiple-hdlss},
for fixed $n$, as $d \rightarrow \infty$, the following results hold.

\begin{enumerate}
\item[(a)]{{If there exists a constant $h$, $1\leq h \leq m$, such that
$\frac{d}{\lambda_h}\rightarrow 0$ and $\frac{d}{\lambda_{h+1}}\rightarrow \infty$,
then $\frac{\hat{\lambda}_j}{\lambda_j }\stackrel{\rm a.s}{\rightarrow}
\frac{\widetilde{Z}_j^T \widetilde{Z}_j}{n}$
for $1\leq j\leq h$,
 where $\widetilde{Z}_j$ is defined in~\eqref{dual:Z},
 and the other $\hat{\lambda}_j$'s satisfy $\frac{\hat{\lambda}_j}{d} \stackrel{\rm a.s}{\rightarrow} K$. In addition,
$\hat{u}_j$ are consistent with rate
$\left( a_j\vee\frac{d}{\lambda_j}\right)^{\frac{1}{2}}$ for $1\leq j\leq h$,
 and the other $\hat{u}_j$'s
are strongly inconsistent with rate
$\left(\frac{\lambda_j}{d}\right)^{\frac{1}{2}}$.}}

\item[(b)] If $\frac{d}{\lambda_1}\rightarrow \infty$,
then the non-zero $\frac{\hat{\lambda}_j}{d}\stackrel{\rm a.s}{\rightarrow} K$, and
the corresponding $\hat{u}_j$ are strongly inconsistent with
 rate $(\frac{\lambda_j}{d})^{\frac{1}{2}}$.

\end{enumerate}
\label{Th:02:hdlss}
\end{theorem}

Some comments are made for the results of Theorem~\ref{Th:02:hdlss}
\begin{itemize}
\item {{If $m=1$, Theorem~\ref{Th:02:hdlss} becomes Theorem~\ref{Th:01:hdlss}.}}

\item {{Even if the non-spike eigenvalues $\lambda_j$, $j>m$, decay slowly,
 the condition $\lambda_1\gg \cdots \gg \lambda_m \gg \lambda_{m+1}$ is enough to guarantee the same properties for $\hat{\lambda}_j$ and
 $\hat{u}_j$, where $1\leq j\leq h$, in Scenario (a).}}

\item  {{If Assumption~\ref{gauss-assumption}  is strengthened to a normality assumption, then
$\frac{\hat{\lambda}_j}{\lambda_j}  \stackrel{\rm a.s}{\rightarrow} \frac{\chi_n^2}{n}$ for $1\leq j\leq h$
 in Scenario (a).}}


\end{itemize}

\subsection{Multiple component spike models with tiered eigenvalues}
\label{ssec:23}

We now consider models where the $m$ eigenvalues can be grouped into $r$ tiers, where the eigenvalues within the same tier are either the same or have the same limit or are of the same order, and the eigenvalues within different tiers have either different limits or are of different orders.

\subsubsection{Cases with increasing sample size $n$} To fix ideas, the first $m$ eigenvalues are grouped into $r$ tiers where there are $q_l(>0)$ eigenvalues in the $l$th tier with $\sum_{l=1}^r q_l=m$. Define $q_0=0$, $q_{r+1}=d-\sum_{l=1}^r q_l$, and the index set of the eigenvalues in the $l$th tier as
\begin{equation}
H_l=\left\{\sum_{k=0}^{l-1} q_k+1,\sum_{k=0}^{l-1}
q_k+2,\cdots, \sum_{k=0}^{l-1} q_k+q_l\right\}, \quad l=1,\cdots,r+1.
\label{index}
\end{equation}
Assume the eigenvalues in the $l$th tier have the same limit $\delta_l(>0)$, i.e. 
\begin{assumption}
$\mbox{\rm lim}_{n\rightarrow \infty} \frac{\lambda_j}{\delta_l}=1, \; j\in H_l, l=1,\cdots,r.$
\label{eigen_assumption}
\end{assumption}
{{The above assumption suggests that it is impossible to separate the sample eigenvectors whose indexes
are in the same tier, and motives us to consider subspace consistency. In addition, we assume that the population eigenvalues from different tiers are
asymptotically different and dominate the other population eigenvalues that
are asymptotically equivalent:}}
\begin{assumption}
as $n \rightarrow \infty$, $\delta_1 > \cdots > \delta_r> \lambda_{m+1}\rightarrow \cdots\rightarrow \lambda_d\rightarrow c_r.$
\label{eigenvalue-assumption-multiple-tiered}
\end{assumption}

{{Under the above setup, we have the following Theorem~\ref{Th:03} which suggests that the eigenvalues with the same limit can not be consistently estimated individually; the corresponding eigenvector estimates are either subspace consistent with the linear space spanned by the eigenvectors, or strongly inconsistent. Similar to the earlier theorems, Theorem~\ref{Th:03} is arranged
according to a decreasing amount of positive information:
\begin{itemize}
\item Theorem~\ref{Th:03}(a): {{If the amount of positive information dominates the amount of negative information up to the $r$th tier, i.e. $\frac{d}{n\delta_r}\rightarrow 0$,
then the estimates for the eigenvectors in the first $r$ tiers  are subspace consistent, and the estimates for the rest are also subspace consistent (but) at a different rate;}}
\item Theorem~\ref{Th:03}(b): {{Otherwise, if the amount of positive information dominates the amount of negative information only up to the $h$th tier $(h\in[1,r])$, i.e. $\frac{d}{n\delta_h}\rightarrow 0$ and $\frac{d}{n\delta_{h+1}}\rightarrow \infty$, then the estimates for the eigenvectors in the first $h$ tiers  are subspace consistent, and the estimates for the rest eigenvectors are strongly-inconsistent;}}
\item Theorem~\ref{Th:03}(c): Finally, if the amount of negative information always dominates, i.e. $\frac{d}{n\lambda_1}\rightarrow \infty$, then the sample eigenvalues are asymptotically indistinguishable, and the sample eigenvectors are strongly inconsistent.
\end{itemize}
}}

In this setting, one key to distinguishing the cases in the theorem is
\begin{equation}
a_l=\mbox{max}_{1\leq k \leq l}\frac{\delta_{k+1}}{\delta_k}, \quad l=1,\cdots,r,
\label{al}
\end{equation}
{{where $\delta_{r+1}=1$, which measures the separation between the sample eigenvectors in the $l$-th tier and those in the first $l-1$ tiers.}} Define the subspace $S_l=\mbox{span}\{u_k, k \in H_l\}$ for $l=1,\cdots,r+1$.

\begin{theorem}  Under Assumptions~\ref{gauss-assumption},~\ref{eigen_assumption} and~\ref{eigenvalue-assumption-multiple-tiered}, as $n \rightarrow \infty$, the following results hold.

\begin{enumerate}
\item[(a)] {{If $\frac{d}{n \delta_r}\rightarrow 0$,
 then $\frac{\hat{\lambda}_j}{\lambda_j} \xrightarrow{\rm a.s} 1$
 for $1\leq j\leq m$.
 In addition, $\hat{u}_j$ are subspace consistent with  $S_l$ 
$j\in H_l, l=1,\cdots, r+1$.}}

\item[(b)] If there exists a constant $h$, $1\leq h \leq r$, such that
$\frac{d}{n\delta_h}\rightarrow 0$ and $\frac{d}{n\delta_{h+1}}\rightarrow \infty$,
then  $\frac{\hat{\lambda}_j}{\lambda_j} \xrightarrow{\rm a.s} 1$ for $ j \in H_l, l=1,\cdots, h$,
and the other non-zero $\hat{\lambda}_j \stackrel{\rm a.s}{\rightarrow} c_\lambda \frac{d}{n} $.
In addition, $\hat{u}_j$ are subspace consistent with $S_l$
 with rate $\left(a_l\vee \frac{d}{n\delta_l}\right)^{\frac{1}{2}}$ for
$j\in H_l, l=1,\cdots, h$, and the other $\hat{u}_j$ are strongly inconsistent
with rate $\left(\frac{n\lambda_j}{d}\right)^{\frac{1}{2}}$.

\item[(c)] If $\frac{d}{n\delta_1}\rightarrow \infty$,
 then the non-zero $\hat{\lambda}_j \stackrel{\rm a.s}{\rightarrow} c_\lambda \frac{d}{n} $,
 and the corresponding $\hat{u}_j$ are strongly inconsistent with rate
$(\frac{n\lambda_j}{d})^{\frac{1}{2}}$.

\end{enumerate}
\label{Th:03}
\end{theorem}

The following comments can be made for the results of Theorem~\ref{Th:03}.
\begin{itemize}

\item If each tier only contains one eigenvalue, i.e. $q_1=\cdots=q_r=1$, then Theorem~\ref{Th:03} simplifies to Theorem~\ref{Th:02}.

\item {{There are additional eigenvalue properties, which are entirely parallel to those in the 2nd remark following Theorem~\ref{Th:02}.
  The corresponding convergence rates in Scenario (a) of Theorem~\ref{Th:03} can be attained
    by replacing $\lambda_j$ by $\delta_j$ in Scenario (a) of  Theorem~\ref{Th:02}.}}

\item {{Assumption~\ref{eigenvalue-assumption-multiple-tiered} can be replaced by
$\delta_1>\cdots >\delta_r \gg \lambda_{m+1}\sim \cdots\lambda_d\sim 1$. Then, the consistency results of the first $r$ tiers of sample eigenvalues in Scenario  (a) or the first $h$ tiers in Scenario (b) remain the same, while
all other results of the form $``\xrightarrow{\rm a.s}"$ for the sample eigenvalues should be replaced by $``\stackrel{\rm a.s}{\sim}"$.
The results for the sample eigenvectors remain same.}}

\item {{Even if the non-spike eigenvalues $\lambda_j$, $j>m$, decay slowly,
 the condition $\delta_1>\cdots >\delta_r \gg \lambda_{m+1}$ is enough to generate the same properties for $\hat{\lambda}_j$ and
  $\hat{u}_j$, with $j\in H_l$, $l\leq r$ as in Scenario (a) and $j\in H_l$, $l\leq h$ as in Scenario (b).}}

\item The cases covered by Theorem~\ref{Th:03} were not studied by Paul (2007)~\cite{paul2007asymptotics}, which required the eigenvalues to be individually estimable.


\item In Theorem~\ref{Th:03}, the dimension $d$ can be fixed. In addition, suppose
$\infty >\delta_1 > \cdots > \delta_r>  \lambda_{m+1}\rightarrow\cdots\rightarrow \lambda_d \rightarrow c_\lambda$ and the eigenvalues satisfying \eqref{eigen_assumption}. Then, the results of Theorem~\ref{Th:03}(a) are consistent with
 the classical asymptotic subspace consistency results implied
 by Theorem 1 of Anderson (1963)~\cite{anderson1963asymptotic}.

\end{itemize}

\subsubsection{Cases with fixed $n$} 
Similar results can be obtained for the fixed $n$ cases (i.e. the HDLSS domain) as summarized below in Theorem~\ref{Th:03:hdlss}. For that, {{we assume that as $d \rightarrow \infty$, the first $m$ eigenvalues fall into $r$ tiers, where the eigenvalues in the same tier are asymptotically equivalent, as stated in the following assumption:}}
\begin{assumption}
$\lambda_j \sim \delta_l, \; j\in H_l, l=1,\cdots,r.$
\label{eigen_assumption2}
\end{assumption}
Different from Assumption~\ref{eigen_assumption} for diverging sample size $n$, now with a fixed $n$, the eigenvalues within the same tier are assumed to be of the same order, rather than of the same limit when $n$ increases to $\infty$. As we will see below in Theorem~\ref{Th:03:hdlss}, one can not separately estimate the eigenvalues of the same order when $n$ is fixed, which is feasible with an increasing $n$ as long as they do not have the same limit as previously shown in Theorem~\ref{Th:03}.

{{In addition, we assume that the population eigenvalues from different tiers are of different orders and dominate the
rest eigenvalues which are asymptotically equivalent:}}
\begin{assumption}
as $d \rightarrow \infty$,  $\delta_1 \gg \cdots \gg \delta_r\gg \lambda_{m+1}\sim\cdots\sim \lambda_d\sim 1.$
\label{eigenvalue-assumption-multiple-tiered-hdlss}
\end{assumption}
Note that for fixed $n$ and $d\rightarrow \infty$, the assumption $\delta_l>\delta_{l+1}$ can not
guarantee asymptotic separation of the corresponding sample
eigenvalues $\hat{\lambda}_j$ for $j\in H_l$ and $\hat{\lambda}_{j}$ for $j\in H_{l+1}$.
Thus, we need to replace Assumption~\ref{eigenvalue-assumption-multiple-tiered}
with Assumption~\ref{eigenvalue-assumption-multiple-tiered-hdlss} in order to asymptotically separate
{{the first $r$ subgroups of sample eigenvalues.
Define
\begin{equation*}
K=\lim_{d\rightarrow \infty} \frac{\sum_{j=m+1}^d \lambda_j}{nd} \quad
{\rm and} \quad A^*_l=\frac{1}{n}\sum_{k\in H_l } \widetilde{Z}_k \widetilde{Z}^T_k, \quad l=1, \cdots, r,
\end{equation*}
which are used to describe the asymptotic properties of the sample eigenvalues in HDLSS settings.}}

\begin{theorem} Under Assumptions~\ref{gauss-assumption},~\ref{eigen_assumption2}
and~\ref{eigenvalue-assumption-multiple-tiered-hdlss}, for fixed $n$, as $d \rightarrow \infty$,
the following results hold.

\begin{enumerate}
\item[(a)] {{If there exists a constant $h$, $1\leq h \leq r$, such that
$\frac{d}{\delta_h}\rightarrow 0$ and $\frac{d}{\delta_{h+1}}\rightarrow \infty$,
then for $ j \in H_l, l=1,\cdots, h$, we have almost surely that
\begin{equation}~\label{A:eigenvalue:hdlss}
 \lambda_{\rm min}(A^*_l) \times {\rm min}_{k\in H_l} \lambda_k \leq  \hat{\lambda}_j  \leq
 \lambda_{\rm max}(A^*_l)\times {\rm max
 }_{k\in H_l} \lambda_k,
\end{equation}
and the other $\hat\lambda_j$'s satisfy $\frac{\hat{\lambda}_j}{d} \stackrel{\rm a.s}{\rightarrow} K$.
In addition, $\hat{u}_j$ are subspace consistent with $S_l$
 with rate $\left(a_l\vee \frac{d}{\delta_l}\right)^{\frac{1}{2}}$ for
$j\in H_l, l=1,\cdots, h$, and the other $\hat{u}_j$'s are strongly inconsistent with
 rate $\left(\frac{\lambda_j}{d}\right)^{\frac{1}{2}}$.}}

\item[(b)] If $\frac{d}{\delta_1}\rightarrow \infty$, then the non-zero
$\frac{\hat{\lambda}_j}{d} \stackrel{\rm a.s}{\rightarrow} K$,
 and the corresponding $\hat{u}_j$ are strongly inconsistent with
rate $\left(\frac{\lambda_j}{d}\right)^{\frac{1}{2}}$.

\end{enumerate}
\label{Th:03:hdlss}
\end{theorem}

The following comments can be made about the results of Theorem~\ref{Th:03:hdlss}.
\begin{itemize}

\item {{If each tier only contains one eigenvalue, i.e. $q_1=\cdots=q_r=1$, then
\eqref{A:eigenvalue:hdlss} becomes $\frac{\hat{\lambda}_j}{\lambda_j}\stackrel{\rm a.s}{\rightarrow}\frac{\widetilde{Z}^T_j\widetilde{Z}_j}{n}$ and
Theorem~\ref{Th:03} becomes Theorem~\ref{Th:02}.
}}
\item {{Even if the non-spike eigenvalues $\lambda_j$, $j>m$, decay slowly,
 the condition $\delta_1\gg\cdots \gg\delta_r \gg \lambda_{m+1}$ can still guarantee the same properties for $\hat{\lambda}_j$ and
  $\hat{u}_j$, with $j\in H_l$, $l\leq h$,  in Scenario (a).}}

\end{itemize}

\section{Discussion}\label{dicussion}

Throughout the paper, we assume that the small eigenvalues have the same limit or the same order as 1, i.e.
$\lambda_{m+1}\rightarrow \cdots \rightarrow \lambda_d \rightarrow c_\lambda$ or
$\lambda_{m+1}\sim \cdots \sim \lambda_d \sim 1$. In fact, this is a convenient WLOG choice. Our results remain valid when these small eigenvalues are not of the same order, and even when some of them are 0. For example, suppose $\lambda_{d_1+1}=\cdots=\lambda_d=0$ for
 $m+1 < d_1 < d$. As shown in Section C of the supplementary material~\citep{shen2011online}, the asymptotic properties of PCA are independent
 of the basis choice for the $d$-dimensional space. If the population eigenvectors $u_j$, $j=1,\ldots,d$, are chosen as the basis of the $d$-dimensional space, the population covariance matrix becomes
  \begin{equation*}
\Sigma=\Lambda=\begin{pmatrix}
  \Lambda_1 &  0_{d_1 \times (d-d_1)} \\
    0_{(d-d_1) \times d_1 }  &  0_{(d-d_1) \times(d-d_1)}
 \end{pmatrix},
 \;\; \mbox{where} \;\; \Lambda_1=\begin{pmatrix}
  \lambda_1 &  \cdots & 0 \\
   \vdots    & \ddots & \vdots  \\
  0  &  \cdots & \lambda_{d_1}
 \end{pmatrix},
 \label{covariancematrix}
\end{equation*}
and $0_{k \times l}$ is the $k$-by-$l$ zero matrix.
Then, the asymptotic properties of PCA under the population covariance matrix $\Sigma$ is the
same as those under the covariance matrix $\Lambda_1$. Therefore, we only need to replace the dimension $d$ by the
{\it effective} dimension $d_1$, and all the earlier results can be obtained.

It would be interesting but challenging to explore the non-asymptotic results
such as large deviations of the angle between the sample  and population eigenvectors.
The properties of sample eigenvectors heavily depend on the sample eigenvalues' properties.
Since we are not aware of any non-asymptotic results for the eigenvalues of the random matrix, then
 it appears to be challenging to obtain non-asymptotic results for sample eigenvectors.

\section{Proofs}\label{proofs}
We now provide detailed proofs for the general Theorem~\ref{Th:03}. To save space, proofs for Theorems~\ref{Th:01},~\ref{Th:01:hdlss},
~\ref{Th:02},~\ref{Th:02:hdlss}, and~\ref{Th:03:hdlss}
 (which are often similar, and simpler) are provided in the supplement~\cite{shen2011online}. We first provide some overview in Section~\ref{sec:overview} and list four lemmas in Section~\ref{proofs1}, and then prove the asymptotic properties of
the sample eigenvalues and the sample eigenvectors in Sections~\ref{proofs2} and~\ref{proofs3}, respectively. 


In this paper, we study the consistency and strong inconsistency of PCA through the angle or the inner product
between a sample eigenvector and the corresponding population eigenvector. We first note that this angle has a nice invariance property: it doesn't
depend on the specific choice of the basis for the $d$-dimensional space, as discussed in details in the supplement~\cite{shen2011online}. Given this invariance property, for the rest of the paper, we choose to use the population eigenvectors $u_j$, $j=1,\ldots,d$, as the basis of the $d$-dimensional {{space, which is equivalent to assuming that $X_i$, $i=1,\ldots,n$, is a $d$-dimensional random vector with mean zero and
a diagonal covariance matrix as $\Sigma=\Lambda={\rm diag}\{\lambda_1, \ldots, \lambda_d\}$. This will simplify our mathematical analysis, see for example~\eqref{interproduct_u} and~\eqref{subspace_u}.}}

We consider general cases where the first $m$ eigenvalues are grouped into $r$ tiers, and WLOG we assume that
$\lambda_1=\cdots=\lambda_{q_1}=\delta_1$, $\cdots$,
$\lambda_{\sum_{l=0}^{r-1}q_l+1}=\cdots=\lambda_{m}=\delta_r$ where $q_0=0$ and $q_l$ are positive integers for $l\ge 1$. In addition, we assume
that each ratio ${\delta_j}/{\delta_i}$, where $ 1\leq i<j\leq r$, converges to a constant less than 1 as $n \rightarrow \infty$. (The following arguments can be extended to cases where only the upper limits of the ratios exist as stated in the theorems, through taking a converging subsequence of the diverging sequence of $n$.)

\subsection{Overview}\label{sec:overview} Our proof makes use of the connection between the sample covariance matrix $\hat{\Sigma}$ and its dual matrix $\hat{\Sigma}_D$, which share the same nonzero eigenvalues.
Since  $\Sigma=\Lambda={\rm diag}\{\lambda_1, \ldots, \lambda_d\}$, then it follows from~\eqref{subguassion:assumption}
and~\eqref{dual:Z} that
the dual matrix can be expressed as
 \begin{equation*}
\hat{\Sigma}_D=n^{-1}X^{T}X=\frac{1}{n}\sum_{j=1}^d \lambda_j\widetilde{Z}_j \widetilde{Z}^T_j,
\end{equation*}
which can be
 rewritten as the sum of two matrices as follows:
\begin{equation}\label{AB}
\hat\Sigma_D=A+B,\quad \mbox{with}\quad A=\frac{1}{n}\sum_{j=1}^m \lambda_j\widetilde{Z}_j \widetilde{Z}^T_j, \quad B=\frac{1}{n}\sum_{j=m+1}^d \lambda_j \widetilde{Z}_j \widetilde{Z}^T_j.
\end{equation}

The proof involves the following several steps. First, we study the asymptotic properties of the eigenvalues of $A$ and $B$ in Lemmas~\ref{lemma:03} and~\ref{lemma:04}, respectively. Then, the \emph{Wielandt's Inequality} (Rao~\cite{rao2002linear}), now restated as Lemma~\ref{lemma:02}, enables us to establish the asymptotic properties of the eigenvalues
of the dual matrix in Section~\ref{proofs2}. Finally, we derive the asymptotic properties of
the sample eigenvectors of $\hat\Sigma$ in Section~\ref{proofs3}. Some intuitive ideas are provided in the supplement~\cite{shen2011online} to help understanding the proof.

\subsection{Lemmas}\label{proofs1}
We list four lemmas that are used in our proof. Lemmas~\ref{lemma:03} and~\ref{lemma:04} are proven in our online supplement, the proofs of which need the following Lemma~\ref{lemma:01} that studies asymptotic properties of the largest and smallest non-zero eigenvalues of
random matrix.

\begin{lemma}
As $n\rightarrow \infty$, the eigenvalues of the matrix $A$ in~\eqref{AB} satisfy
\begin{equation*}
\frac{\lambda_j(A)}{\lambda_j}\xrightarrow{\rm a.s} 1, \quad \mbox{for}\quad  j=1,\cdots, m,
\end{equation*}
where $\lambda_j(A)$ denotes the $j$th largest eigenvalue of the matrix $A$.
\label{lemma:03}
\end{lemma}

\begin{lemma}
{{As $n\rightarrow \infty$, the eigenvalues of the matrix $B$ in~\eqref{AB} satisfy that,
for $j=1,\cdots, [n\wedge(d-m)]$,}}
\begin{eqnarray}\label{B:n/d8}
\frac{\lambda_j(B)}{\lambda_{j+m}}\stackrel{\rm a.s}{\rightarrow}  1,  \quad   {\rm for} \quad \frac{d}{n}\rightarrow 0,\\\label{B:n/d0}
\frac{\lambda_j(B)}{\lambda_{j+m}}\stackrel{\rm a.s}{\rightarrow}  \frac{d}{n},  \quad  {\rm for} \quad \frac{d}{n}\rightarrow \infty,
\end{eqnarray}
and almost surely,
\begin{equation}
\label{B:n/dother}
\overline{{\rm lim}}\lambda_1(B) \leq  c \times \overline{{\rm lim}}\frac{d}{n},  \quad otherwise,
\end{equation}
{{where $c$ is a constant.}}
\label{lemma:04}
\end{lemma}

\begin{remark}
{{If $\lambda_{m+1}\rightarrow \cdots \rightarrow \lambda_d$ is relaxed to $\lambda_{m+1} \sim \cdots \sim \lambda_d$,
then  $``\stackrel{\rm a.s}{\rightarrow}"$ is replaced by
  $``\stackrel{\rm a.s}{\sim}"$ in~\eqref{B:n/d8} and~\eqref{B:n/d0}.}}
\end{remark}

\begin{lemma}~\label{lemma:01}
{{Suppose $B=\frac{1}{q}VV^T$ where $V$ is an $p\times q$ random matrix composed
of i.i.d. random variables with zero mean, unit variance and finite fourth moment.
As $q\rightarrow \infty$ and $\frac{p}{q}\rightarrow c \in [0, \infty)$, the largest and smallest
non-zero eigenvalues of $B$ 
converge almost surely to $(1+\sqrt{c})^2$ and $(1-\sqrt{c})^2$, respectively.}}
\end{lemma}

\begin{remark}
{{Lemma~\ref{lemma:01} is known as the Bai-Yin's law~\cite{bai1993limit}.
As in Remak 1 of~\cite{bai1993limit}, the smallest non-zero eigenvalue is
the $p-q+1$ smallest eigenvalue of $B$ for $c>1$.}}
\end{remark}

\begin{lemma} (Wielandt's Inequality~\cite{rao2002linear}).
If $A, B$ are $p\times p $ real symmetric matrices, then for all $j = 1,\ldots,p$,
\begin{equation*}
\left\{ \begin{array}{lll}
  \lambda_j(A)&+& \lambda_p(B)\\
   \lambda_{j+1}(A)&+& \lambda_{p-1}(B)\\
     &\vdots& \\
   \lambda_p(A)&+& \lambda_j(B)
 \end{array} \right\}
 \leq \lambda_j(A+B) \leq
 \left\{ \begin{array}{lll}
  \lambda_j(A)&+& \lambda_1(B)\\
   \lambda_{j-1}(A)&+& \lambda_2(B)\\
     &\vdots& \\
   \lambda_1(A)&+& \lambda_j(B)
 \end{array}\right\}.
\end{equation*}
\label{lemma:02}
\end{lemma}

\subsection{Asymptotic properties of the sample eigenvalues}\label{proofs2}
We now study the asymptotic properties of the sample eigenvalues $\hat{\lambda}_j$, for $j=1,\cdots, [n\wedge d]$, which are the same as the eigenvalues of the dual matrix $\hat{\Sigma}_D$, denoted as $\lambda_j(\hat{\Sigma}_D)=\lambda_j(A+B)$.

\subsubsection{Scenario (a) in Theorem~\ref{Th:03}}\label{Scenario (a):eigenvalue}

{{Note that $\frac{d}{n\lambda_m}\rightarrow 0$ ($\delta_r=\lambda_m$) contains three different cases:
$\frac{d}{n}\rightarrow 0$, $\infty$ or $\overline{{\rm lim}}\frac{d}{n}=c^*$ ($0<c^*\leq \infty$).
The proofs are different for each case and are provided separately below.}}


{{Consider the first one: $\frac{d}{n}\rightarrow 0$. If in addition we have $\lambda_m\rightarrow \infty$,
then Lemma~\ref{lemma:02} suggests that
\begin{equation}~\label{A+B:1:m}
\frac{\lambda_j(A)}{\lambda_j} \leq \frac{\hat{\lambda}_j}{\lambda_j}\leq \frac{\lambda_j(A)}{\lambda_j}+\frac{\lambda_1(B)}{\lambda_j},
\end{equation}
which, together with $\lambda_m\rightarrow \infty$, \eqref{B:n/d8} and Lemma~\ref{lemma:03}, yields that
\begin{equation}~\label{engvalue1st:consistency}
\frac{\hat{\lambda}_j}{\lambda_j} \stackrel{\rm a.s}{\rightarrow} 1, \quad j=1,\cdots, m.
\end{equation}

Instead, if $\lambda_m<\infty$, according to Theorem 1 ($c=0$) of~\cite{baik2006eigenvalues}, we still
have~\eqref{engvalue1st:consistency}.
In addition, according to Lemma~\ref{lemma:02}, we have
that
\begin{equation} \label{eigenvalue:A+B:>m}
\frac{\lambda_j(B)}{\lambda_j}
\leq \frac{\hat{\lambda}_j}{\lambda_j}\leq
 \frac{\lambda_j(A)}{\lambda_j}+\frac{\lambda_1(B)}{\lambda_j},
\end{equation}
which, together with \eqref{B:n/d8},  $\lambda_j(A)=0$ for $j\geq m+1$ and $\lambda_{m+1} \rightarrow \lambda_d\rightarrow c_\lambda $
yields that
\begin{eqnarray}~\label{engvalue>m:consistency}
&& \frac{\hat{\lambda}_j}{\lambda_j} \stackrel{\rm a.s}{\rightarrow } 1, \quad j=m+1,\cdots, [n\wedge (d-m)],\\ \nonumber
&& \hat{\lambda}_j ={\rm O}_{\rm a.s} (1), \quad j=[n\wedge (d-m)]+1,\cdots, [n\wedge d].
\end{eqnarray}}}

{{Now, consider the second case: $\frac{d}{n}\rightarrow \infty$. Since $\frac{d}{n\lambda_m}\rightarrow 0$, then $\lambda_m\rightarrow \infty$,
which, together with \eqref{B:n/d0}, \eqref{A+B:1:m} and Lemma~\ref{lemma:03}, yields~\eqref{engvalue1st:consistency}. In addition, it
 follows from~\eqref{B:n/d0},~\eqref{eigenvalue:A+B:>m},
$\lambda_j(A)=0$ for $j\geq m+1$ and $\lambda_{m+1}\rightarrow \lambda_d $ that
\begin{equation}~\label{engvalue>m:inconsistency}
\frac{\hat{\lambda}_j}{\lambda_j} \stackrel{\rm a.s}{\rightarrow} \frac{d}{n}, \quad j=m+1,\cdots, [n\wedge d].
\end{equation}}}

{{Finally, consider the third case: $\overline{{\rm lim}}\frac{d}{n}=c^*$ ($0<c^*\leq \infty$). Similarly, it follows from
 $\frac{d}{n\lambda_m}\rightarrow 0$ that $\lambda_m\rightarrow \infty$, which,
jointly with \eqref{B:n/dother}, \eqref{A+B:1:m} and Lemma~\ref{lemma:03},
yields~\eqref{engvalue1st:consistency}. In addition, note that~\eqref{B:n/dother},~\eqref{eigenvalue:A+B:>m},
$\lambda_j(A)=0$ for $j\geq m+1$ and $\lambda_{m+1}\rightarrow \lambda_d\rightarrow c_\lambda$, and then
almost surely we have
\begin{equation}~\label{engvalue>m:other}
\overline{{\rm lim}}\hat{\lambda}_j \leq  c \times \overline{{\rm lim}}\frac{d}{n}, \quad j=m+1,\cdots, [n\wedge d],
\end{equation}
where $c$ is a constant.}}


{{All together, we have proved the consistency of the first $m$ sample eigenvalues in~\eqref{engvalue1st:consistency}
and the asymptotic properties of the rest of the non-zero sample eigenvalues
in~\eqref{engvalue>m:consistency},
\eqref{engvalue>m:inconsistency} and~\eqref{engvalue>m:other} for Scenario (a).}}

\subsubsection{Scenario (b) in Theorem~\ref{Th:03}}\label{Scenario (b):eigenvalue}

{{Since $\frac{d}{n\delta_{h+1}} \rightarrow \infty$, then $\frac{d}{n} \rightarrow \infty$.
According to~\eqref{B:n/d0}, ~\eqref{A+B:1:m} and $\frac{d}{n\delta_h} \rightarrow 0$, we have
\begin{equation}~\label{engvalue1st:consistency:(b)}
\frac{\hat{\lambda}_j}{\lambda_j} \stackrel{\rm a.s}{\rightarrow} 1, \quad j\in H_l, l=1,\cdots, h.
\end{equation}
In addition, note that
\begin{equation}~\label{eigenvalue:j>h}
\frac{\lambda_j(B)n}{d}
\leq \frac{\hat{\lambda}_j n}{d}\leq
 \frac{\lambda_j(A)n}{d}+\frac{\lambda_1(B) n}{d},
\end{equation}
which, together with $\frac{d}{n\delta_{h+1}} \rightarrow \infty$,~\eqref{B:n/d0} and Lemma~\ref{lemma:03}, yields
that
\begin{equation*}
\hat{\lambda}_j \stackrel{\rm a.s}{\rightarrow} c_\lambda \frac{d}{n}, \quad j=\sum_{k=1}^h q_k+1,\cdots, [n\wedge d].
\end{equation*}}}

\subsubsection{Scenario (c) in Theorem~\ref{Th:03}}\label{Scenario (c):eigenvalue}

{{Since $\frac{d}{n\delta_1} \rightarrow \infty$, then $\frac{d}{n} \rightarrow \infty$.
Then it follows from~\eqref{B:n/d0},~\eqref{eigenvalue:j>h} and $\frac{d}{n\delta_1}\rightarrow \infty$ that
\begin{equation*}
  \hat{\lambda}_j \stackrel{\rm a.s}{\rightarrow} c_\lambda \frac{d}{n}, \quad j=1,\cdots, [n\wedge d].
\end{equation*}}}

\subsection{Asymptotic properties of the sample eigenvectors}\label{proofs3}
We first state two results that simplify the proof. As aforementioned, in light of the invariance property of the angle, we choose the population eigenvectors $u_j$, $j=1,\ldots, d$,
as the basis of the $d$-dimensional space. It then follows that $u_j=e_j$ where the $j$th
component of $e_j$ equals to 1 and all the other components equal to zero. This suggests that
\begin{equation}
\mid<\hat{u}_j,u_j>\mid^2=\mid<\hat{u}_j,e_j>\mid^2=\hat{u}^2_{j,j},
\label{interproduct_u}
\end{equation}
and for any index set $H$,
\begin{eqnarray}
\mbox{cos}\left[\mbox{angle}\left(\hat{u}_j, \mbox{span}\{u_k, k\in H\}\right)\right]=\sum_{k\in H} \hat{u}^2_{k,j}.
\label{subspace_u}
\end{eqnarray}

As a reminder, the population eigenvalues are grouped into $r+1$ tiers and the index set of the eigenvalues in the $l$th tier $H_l$ is defined in \eqref{index}. Define
\begin{equation}
\hat{U}_{k,l}=(\hat{u}_{i,j})_{i\in H_k, j\in H_l}, \quad 1\leq k, l \leq r+1.
\label{Uij}
\end{equation}
Then, the sample eigenvector matrix $\hat{U}$
can be rewritten as the following:
\begin{equation*}
\hat{U}=[\hat{u}_1,\hat{u}_2,\cdots,\hat{u}_d]=\begin{pmatrix}
  \hat{U}_{1,1} & \hat{U}_{1,2} & \cdots & \hat{U}_{1,r+1} \\
    \hat{U}_{2,1} & \hat{U}_{2,2} & \cdots & \hat{U}_{2,r+1} \\
   \vdots    & \vdots & & \vdots  \\
  \hat{U}_{r+1,1} & \hat{U}_{r+1,2} & \cdots & \hat{U}_{r+1,r+1}
 \end{pmatrix}.
\end{equation*}

To derive the asymptotic properties of the sample eigenvectors $\hat{u}_j$, we consider the three scenarios of Theorem~\ref{Th:03} separately.



\subsubsection{Scenario (b) in Theorem~\ref{Th:03}}\label{Scenario (b)}
Under this scenario, there exists a constant $h\in [1,r]$,
such that $\frac{d}{n\delta_h}\rightarrow 0$ and $\frac{d}{n\delta_{h+1}}\rightarrow \infty$.
From~\eqref{subspace_u}, to
show the subspace consistency with $S_l$ and rate $\left(a_l\vee \frac{d}{n\delta_l}\right)^{\frac{1}{2}}$,
we only need to show that
\begin{equation}
\sum_{k\in H_l}\hat{u}^2_{k,j}=1+{\rm o}_{a.s}(a_l)\vee{\rm O}_{a.s}(\frac{d}{n\delta_l}) , \quad j\in H_l, l=1,\cdots, h,
\label{eigenvector:H1}
\end{equation}
where, as defined in~\eqref{al} in Section~\ref{ssec:23}, $a_l=\mbox{max}_{1\leq k \leq l}\frac{\delta_{k+1}}{\delta_k}$, $l=1,\cdots, r$.
Below we provide the proof for $l=1$. The process is similar for $l=2,\cdots, h$, which is omitted to save space.

Note that for $l=1$,
the left hand side of~\eqref{eigenvector:H1} becomes the sum of squares of the column elements
in the matrix $\hat{U}_{1,1}$ (defined in~\eqref{Uij}).
Thus, to prove \eqref{eigenvector:H1}, 
we first show that this sum of squares converges to 1, and then
establish the convergence rate $a_1\vee \frac{d}{n\delta_1}$.

For the first step, let $Z=(Z_1,\cdots,Z_n)$, where $Z_i=(z_{i,1},\cdots, z_{i,d})^T=\Lambda^{-\frac{1}{2}}X_i$
from~\eqref{subguassion:assumption}.
Denote $S=\Lambda^{-\frac{1}{2}}\hat{U}\hat{\Lambda}^{\frac{1}{2}}$ where
$\hat{U}$ is the sample eigenvector matrix and $\hat{\Lambda}$ is the sample eigenvalue matrix
defined in~\eqref{sample:covariance}. 
We can show that
$SS^T=\frac{1}{n}ZZ^T.$
Considering the $k$-th diagonal entry of the matrices on the two sides and noting that
$s_{k,j}=\lambda^{-\frac{1}{2}}_k \hat{\lambda}^{\frac{1}{2}}_j \hat{u}_{k,j}$,
we have the following
\begin{equation}
\frac{1}{n}\sum_{i=1}^n z^2_{i,k}=\sum_{j=1}^{d} s^2_{k,j} =
\lambda^{-1}_k
\sum_{j=1}^{d} \hat{\lambda}_j\hat{u}^2_{k,j}, \quad k=1,\cdots,d.
\label{diagk}
\end{equation}

{{Select the first $[n\wedge d]$ rows of $Z$ and denote the resulting random matrix as $Z^*$. Then, we have ${\rm max}_{1\leq k\leq [n\wedge d]} \frac{1}{n}\sum_{i=1}^n z^2_{i,k}
\leq \lambda_{\rm max}(\frac{1}{n}Z^* {Z^*}^T )$. Note that $\frac{d}{n}\rightarrow \infty$ here, so
$[n\wedge d]=n$. According to Lemma~\ref{lemma:01}, we have $\lambda_{\rm max}(\frac{1}{n}Z^* {Z^*}^T )
\stackrel{\rm a.s}{\rightarrow} 4$, which suggests that almost surely $\hat{u}^2_{k,j}
\leq 4 \frac{\lambda_k}{\hat\lambda_j}$ for $j=1,\cdots, [n\wedge d]$,
as $n\rightarrow \infty$.  Then, 
given the
asymptotic properties of $\hat\lambda_j$ in Scenario (b) of Theorem~\ref{Th:03} (Section~\ref{proofs2}),
it follows that}}
\begin{equation}
\hat{u}^2_{k,j}=
\left\{ \begin{array}{ll}
   {\rm O}_{\rm a.s}(\frac{\lambda_k}{\lambda_j}) &  \quad j\in H_l, l=1,\cdots,h,\\
 {\rm O}_{\rm a.s}(\frac{n\lambda_k}{d}) &   \quad j=\sum_{l=1}^h q_l+1,\cdots,[n\wedge d].
 \end{array} \right.
 \label{estimate:h}
\end{equation}

In addition, the $k$th diagonal entry of $S^T S$ is less than or equal to its largest eigenvalue,
i.e. the largest eigenvalue of $\frac{1}{n}Z^TZ$.
Hence, we have
\begin{equation}
\hat{\lambda}_j
\sum_{k=1}^{d} \lambda^{-1}_k\hat{u}^2_{k,j}=\sum_{k=1}^{d} s^2_{k,j}\leq
\lambda_{\mbox{max}}(\frac{1}{n}Z^TZ), \quad j \in H_l, l=1,\cdots,h.
\label{diagkk}
\end{equation}
 According to Lemma~\ref{lemma:01} and $\frac{d}{n}\rightarrow \infty$,
 we have that
 \begin{equation}
 \lambda_{\mbox{max}}(\frac{1}{n}ZZ^T)\stackrel{\rm a.s}{\rightarrow}(\frac{d}{n}).
 \label{eigen:max}
 \end{equation}
From~\eqref{engvalue1st:consistency:(b)}, \eqref{diagkk}, \eqref{eigen:max} and  $\lambda_{m+1}\rightarrow \lambda_d \rightarrow c_\lambda$, we have
\begin{equation}
\sum_{k=m+1}^{d}\hat{u}^2_{k,j}={\rm O}_{\rm a.s}(\frac{d}{n\lambda_j}), \quad j\in H_l, l=1,\cdots,h.
\label{estimate:m+1}
\end{equation}
Note that $\lambda_j\ll\frac{d}{n}$, for $j=\sum_{l=1}^h q_l+1,\cdots, m$, which together with \eqref{estimate:h}
and \eqref{estimate:m+1}, yields that
\begin{equation}
\sum_{k=\sum_{t=1}^h q_t+1}^{d}\hat{u}^2_{k,j}={\rm O}_{\rm a.s}
\left(\frac{d}{n\lambda_j}\right), \quad j \in H_l, l=1,\cdots,h.
\label{estimate:h+1}
\end{equation}

According to~\eqref{diagk}  and $\lambda_k=\delta_1$, $k\in H_1$,
 we obtain that for $k\in H_1$,
\begin{eqnarray}
\nonumber
 \frac{1}{n}\sum_{i=1}^n z^2_{i,k}=
\lambda^{-1}_k\sum_{j=1}^d \hat{\lambda}_j\hat{u}^2_{k,j}
%
&\leq & \delta^{-1}_1\sum_{j\in H_1}\hat{\lambda}_1\hat{u}^2_{k,j}+
\delta^{-1}_1\sum_{j \notin H_1}\hat{\lambda}_{q_1+1}\hat{u}^2_{k,j}\\
\label{du11}
&=&\delta^{-1}_1(\hat{\lambda}_1-\hat{\lambda}_{q_1+1})\sum_{j\in H_1}\hat{u}^2_{k,j}+
\delta^{-1}_1\hat{\lambda}_{q_1+1}.
\end{eqnarray}
In addition, it follows from \eqref{engvalue1st:consistency:(b)}
that  $\delta^{-1}_1(\hat{\lambda}_1-\hat{\lambda}_{q_1+1})\xrightarrow{\rm a.s}
\delta^{-1}_1(\delta_1-\delta_2)=(1-c)$, and
$\delta^{-1}_1\hat{\lambda}_{q_1+1}\xrightarrow{\rm a.s} c$,
where $c=\mbox{lim}_{n\rightarrow \infty} \frac{\delta_2}{\delta_1}<1$.

Note that $\frac{1}{n}\sum_{i=1}^n z^2_{i,k}=1+{\rm o}_{a.s}(1)$, which together with~\eqref{du11},
yields that
\begin{eqnarray*}
\nonumber
1+{\rm o}_{a.s}(1)&\leq& (1-c)\underline{\mbox{lim}}_{n \rightarrow \infty}\sum_{j\in H_1}\hat{u}^2_{k,j}+c\\
&\leq&(1-c)\overline{\mbox{lim}}_{n \rightarrow \infty}\sum_{j\in H_1}\hat{u}^2_{k,j}+c \leq 1,
\end{eqnarray*}
which yields
$\sum_{j\in H_1}\hat{u}^2_{k,j} \xrightarrow{\rm a.s}  1$, $k\in H_1$.
The above means that the sum of squares of the row elements of $\hat{U}_{1,1}$ converges to 1.
Given that the sample eigenvectors all have norm 1, the sum of squares of the row or the column elements of $\hat{U}_{1,1}$ is less than
or equal to 1. It then follows that the sum of squares of the column elements of $\hat{U}_{1,1}$
converges to 1, which finishes the first step of the proof.

For the second step of the proof, we need to establish
the convergence rate $a_1\vee \frac{d}{n\delta_1}$ of the above sum of squares.
Having shown that the sum of squares of the row elements of
$\hat{U}_{1,1}$ converges to 1, it follows that
the sum of squares of the row elements of $\hat{U}_{1, 2}$ converges to
0. Furthermore, the sum of the squares of the column elements of
$\hat{U}_{1, 2}$ converges to 0, as follows:
\begin{equation}
\sum_{k\in H_1}\hat{u}^2_{k,j}= {\rm o}_{a.s}(1), \quad j\in H_2.
\label{vanish1}
\end{equation}
WLOG, we assume that $\frac{\delta_3}{\delta_2}\rightarrow 0$.
(If the limit is greater than 0, we can combine the index sets $H_2$ and $H_3$ together
to check whether $\frac{\delta_4}{\delta_2}\rightarrow 0$ converges to 0.
If not, we keep combining the index sets together until the big jump appears.)
Given that $\frac{\delta_3}{\delta_2}\rightarrow 0$, ~\eqref{estimate:h} and~\eqref{estimate:h+1}, it follows
that
\begin{equation}
\sum_{k\in H_3\cup\cdots \cup H_{r+1}}\hat{u}^2_{k,j}= {\rm o}_{a.s}(1), \quad j\in H_2.
\label{vanish3}
\end{equation}
From~\eqref{vanish1} and~\eqref{vanish3}, we have that
\begin{equation*}
\sum_{k\in H_2}\hat{u}^2_{k,j}=1+ {\rm o}_{a.s}(1), \quad j\in H_2,
\label{H2}
\end{equation*}
which means that the sum of squares of the column elements of $\hat{U}_{2,2}$
also converges to 1. Again, since the sum of squares of the row or column elements of
$\hat{U}_{2,2}$ is less than or equal to 1, it follows that
the sum of squares of the row elements of $\hat{U}_{2,2}$ must converge to 1:
\begin{equation}
\sum_{j\in H_2}\hat{u}^2_{k,j}=1+ {\rm o}_{a.s}(1), \quad k\in H_2.
\label{H2r}
\end{equation}

Given that $\hat{\lambda}_j\xrightarrow{a.s} \lambda_j=\delta_2$, $j \in H_2$, and \eqref{H2r}, it follows that, for
$k \in H_2$,
\begin{eqnarray*}
\nonumber
&& 1+{\rm o}_{a.s}(1)=\frac{1}{n}\sum_{i=1}^n z^2_{i,k}=\lambda^{-1}_k\sum_{j=1}^d \hat{\lambda}_j\hat{u}^2_{k,j}\\
&& \geq
\delta^{-1}_2\sum_{j\in H_1} \hat{\lambda}_j\hat{u}^2_{k,j}+ \delta^{-1}_2\sum_{j\in H_2} \hat{\lambda}_j\hat{u}^2_{k,j}=\delta^{-1}_2\sum_{j\in H_1} \hat{\lambda}_j\hat{u}^2_{k,j}+1+{\rm o}_{a.s}(1),
\label{u22}
\end{eqnarray*}
which yields $\delta^{-1}_2\sum_{j \in H_1} \hat{\lambda}_j\hat{u}^2_{k,j}={\rm o}_{a.s}(1)$, $k\in H_2$.

For $j \in H_1$, we have that $\hat{\lambda}_j\xrightarrow{a.s} \lambda_j=\delta_1$; hence, it follows that
\begin{equation*}
\sum_{j\in H_1}\hat{u}^2_{k,j}={\rm o}_{a.s}(\frac{\delta_2}{\delta_1}), \quad k \in H_2,
\end{equation*}
which yields that
\begin{equation}
\sum_{k\in H_2}\hat{u}^2_{k,j}={\rm o}_{a.s}(\frac{\delta_2}{\delta_1}), \quad j \in H_1.
\label{VH_2}
\end{equation}
In addition, from~\eqref{estimate:h} and~\eqref{estimate:h+1}, we have
\begin{equation}
\sum_{k\in H_3\cup\cdots \cup H_{r+1}}\hat{u}^2_{k,j}= {\rm o}_{a.s}(\frac{\delta_2}{\delta_1}), \quad j\in H_1.
\label{H>3}
\end{equation}

From~\eqref{VH_2}, ~\eqref{H>3} and $\frac{\delta_2}{\delta_1}\gg \frac{d}{n\delta_1}$, we have
\begin{equation*}
\sum_{k\in H_1}\hat{u}^2_{k,j}=1+{\rm o}_{a.s}(\frac{\delta_2}{\delta_1})=1+{\rm o}_{a.s}(a_1)\vee{\rm O}_{a.s}(\frac{d}{n\delta_1}), \quad j\in H_1,
\end{equation*}
which suggests that the sum of squares of the column elements
of $\hat{U}_{1,1}$ converges to 1 with
the convergence rate $a_1\vee \frac{d}{n\delta_1}$, as stated in
\eqref{eigenvector:H1} for $l=1$. The proof of \eqref{eigenvector:H1}
is similar for  $l=2,\cdots, h$. Thus, we have shown the subspace consistency portion of the results in
Scenario (b).

Finally, the strong inconsistency in Scenario (b) follows directly from~\eqref{estimate:h} by setting $k=j$:
\begin{equation}\label{eigenvector:inconsistency:rest}
\mid<\hat{u}_j, u_j>\mid^2=\hat{u}^2_{j,j}={\rm O}_{\rm a.s}\left(\frac{n\lambda_j}{d} \right), \quad
j=\sum^h_{l=1} q_l+1, \cdots, [n\wedge d].
\end{equation}
Hence, we have finished the proof of Scenario (b) in Theorem~\ref{Th:03}.


\subsubsection{Scenario (a) in Theorem~\ref{Th:03}}

{{As in Section~\ref{Scenario (a):eigenvalue}, $\frac{d}{n\lambda_m}\rightarrow 0$ ($\delta_r=\lambda_m$) contains three different cases:
$\frac{d}{n}\rightarrow 0$, $\infty$, or $\overline{{\rm lim}}\frac{d}{n}=c^*$ ($0<c^*\leq \infty$), which we shall prove separately.}}


{{Consider that $\frac{d}{n}\rightarrow 0$, then
\eqref{eigenvector:H1} in Section~\ref{Scenario (b)} becomes
\begin{equation}
\sum_{k\in H_l}\hat{u}^2_{k,j}=1+{\rm o}_{a.s}(a_l)\vee{\rm O}_{a.s}(\frac{1}{\delta_l}) , \quad j\in H_l, l=1,\cdots, r,
\label{eigenvector:H1a}
\end{equation}
which, together with similar arguments as those in proving Scenario (b), leads to
\begin{equation}
\sum_{k\in H_{r+1}}\hat{u}^2_{k,j}=1+{\rm o}_{a.s}(a_r)\vee{\rm O}_{a.s}(\frac{1}{\delta_r}) , \quad m+1\leq j \leq [n\wedge d].
\label{eigenvector:H2a}
\end{equation}
Note that if $\delta_l<\infty$, then $``{\rm O}_{a.s}(\frac{1}{\delta_l})"$ in~\eqref{eigenvector:H1a} and~\eqref{eigenvector:H2a} becomes
$``{\rm o}_{a.s}(1)"$.}}

{{For the second case: $\frac{d}{n}\rightarrow \infty$, we have that $`` {\rm O}_{a.s}(\frac{1}{\delta_l}) "$
in~\eqref{eigenvector:H1a} and~\eqref{eigenvector:H2a}
are replaced by  $`` {\rm O}_{a.s}(\frac{d}{n\delta_l}) "$.}}

Finally, consider the third case: $\overline{{\rm lim}}\frac{d}{n}=c^*$ ($0<c^*\leq \infty$). If $c^*<\infty$, then
the convergence rates are the same as those in~\eqref{eigenvector:H1a} and~\eqref{eigenvector:H2a}.
Otherwise, 
$`` {\rm O}_{a.s}(\frac{1}{\delta_l}) "$
in~\eqref{eigenvector:H1a} and~\eqref{eigenvector:H2a}
are replaced by  $`` {\rm O}_{a.s}(\frac{d_n}{n\delta_l}) "$, where the sequence
$\{d_n\}$ is defined in Section E.2 of the Supplement~\cite{shen2011online}.



\subsubsection{Scenario (c) in Theorem~\ref{Th:03}}

Finally, for Scenario (c) where $\frac{d}{n\delta_1}\rightarrow 0$, the strong inconsistency
 in Theorem~\ref{Th:03} follows from~\eqref{estimate:h} by setting $k=j$.


\begin{supplement}
\stitle{Additional Proofs} 
\sdescription{Detailed proofs are provided for Theorems~3.1, 3.2, 3.3, 4.1, 4.2, 4.4, and the necessary lemmas.}
\slink[url]{http://www.unc.edu/~dshen/PCA/PCASupplment.pdf}
\end{supplement}

\bibliographystyle{imsart-number}
\bibliography{PCA}

\end{document}